\documentclass[12pt]{article}
\usepackage{graphicx, amsmath, amssymb, cite, setspace, color, relsize, bm, bbold, array, hyperref, dsfont, bbold, cancel,soul,placeins,hyperref,amsfonts,pifont,comment}
\usepackage{wasysym}
\usepackage[dvipsnames]{xcolor}

\usepackage[top=1 in, bottom=1 in, left=1 in, right=1 in]{geometry}

\def\sigmaTWO{\textrm{Tr}[\hat{\sigma}^2]}
\def\sigmaTHREE{\textrm{Tr}[\hat{\sigma}^3]}
\def\sigmaFOUR{\textrm{Tr}[\hat{\sigma}^4]}

\usepackage{mathtools}
\usepackage{amsthm}
\usepackage{amssymb}
\usepackage[alphabetic]{amsrefs}
\usepackage{graphicx}
\usepackage{tikz}
	\usetikzlibrary{decorations.pathreplacing}
	\usetikzlibrary{patterns}
\usepackage{caption}
\usepackage{enumitem}

\newcommand{\captionfonts}{\footnotesize} 
\makeatletter
\long\def\@makecaption#1#2{%
  \vskip\abovecaptionskip
  \sbox\@tempboxa{{\captionfonts #1: #2}}%
  \ifdim \wd\@tempboxa >\hsize
    {\captionfonts #1: #2\par}
  \else
    \hbox to\hsize{\hfil\box\@tempboxa\hfil}%
  \fi
  \vskip\belowcaptionskip}
\makeatother

\def\lsim{ \lower .75ex \hbox{$\sim$} \llap{\raise .27ex
\hbox{$<$}} }
\def\gsim{ \lower .75ex \hbox{$\sim$} \llap{\raise .27ex
\hbox{$>$}} }

\usepackage{eso-pic}
\usepackage{lipsum}

\let\oldsqrt\sqrt
\def\sqrt{\mathpalette\DHLhksqrt}
\def\DHLhksqrt#1#2{%
\setbox0=\hbox{$#1\oldsqrt{#2\,}$}\dimen0=\ht0
\advance\dimen0-0.2\ht0
\setbox2=\hbox{\vrule height\ht0 depth -\dimen0}%
{\box0\lower0.4pt\box2}}

\begin{document}

\begin{centering}
{\Large {An Improved Long-Time Bishop-Gromov Theorem Using Shear}}\\
\vspace{6mm} 

{\large {Adam R.~Brown$^{a,b}$ and Michael H.~Freedman$^{c,d}$}}\\
\end{centering} 
\vspace{2mm} 
\begin{centering}
$^{a}$ {Google Research (Blueshift), Mountain View, California}

$^{b}$ {Physics Department, Stanford University, Stanford, California }

$^{c}$ {Microsoft Research, Santa Barbara, California}

$^{d}$ {Mathematics Department, UC Santa Barbara, California }

\end{centering}

\begin{abstract}
\noindent The Bishop-Gromov theorem is a comparison theorem of differential geometry that upperbounds the growth of volume of a geodesic ball in a curved space. For many spaces, this bound is far from tight. We identify a major reason the bound fails to be tight:  it neglects the effect of \emph{shear}. By using higher curvature invariants to lowerbound the average shear, we are able to place tighter-than-Bishop-Gromov upperbounds on the late-time growth rates of geodesic balls in homogeneous spaces with non-positive sectional curvature. We also provide concrete guidance on how our theorem can be generalized to inhomogeneous spaces, to spaces with positive sectional curvatures, and to intermediate and short times. In \cite{BG1} we discovered an enhancement to the BG theorem that was strongest at early times, and that relied upon additive properties of families of Jacobi equations; in this paper we find a different enhancement at late times that connects to multiplicative properties of families of Jacobi equations. A novel feature shared by both papers is the consideration of families of equations that are not coupled but whose coefficients are correlated. 
\end{abstract}

\newpage

\tableofcontents 
\section{Introduction}

\subsection{The Bishop-Gromov bound}
A geodesic ball is formed by picking a starting point in a Riemannian manifold and then shooting out geodesics in every direction. The geodesic ball of radius $t$ is the union of all the points that have been touched at least once by geodesics of length  $t$.  We will be interested in the volume of a geodesic ball as a function of its radius. 

For a given space, the volume of a geodesic ball can in principle be calculated exactly. However, this will require solving the geodesic equation, which will typically be prohibitively computationally expensive. For that reason, we will seek simple upperbounds on the volume.

The most famous such upperbound is the Bishop-Gromov (BG) bound. Let's describe the BG bound for the case of a homogeneous space, which means a space for which all points are the same (and for which, therefore, the volume of a geodesic ball does not depend on its starting point).  Step one is to calculate the most negative component of the Ricci tensor by minimizing $\mathcal{R}_{\mu \nu} X^{\mu} X^{\nu}$ over all unit tangent vectors $X^{\mu}$ to give  $\mathcal{R}_{\mu \nu} X_\textrm{min}^{\mu} X_\textrm{min}^{\nu}$. Step two is to consider the maximally symmetric space with the same dimension, scaled so that its Ricci curvature is equal to the minimal Ricci curvature of the original space, $\mathcal{R}_{ab} =  
(\mathcal{R}_{\mu \nu} X_\textrm{min}^{\mu} X_\textrm{min}^{\nu}) g_{ab}$. Finally, we calculate the volume of a geodesic ball in this maximally symmetric space, which is given by
\begin{equation}
\textrm{BG}(t)  \equiv   \Omega_{d-1} \, \int_0^t d\tau \, \operatorname{sn}\Bigl(  \frac{{ \mathcal{R}_{\mu \nu} X_\textrm{min}^{\mu} X_\textrm{min}^{\nu}}}{{d-1}}, \tau \Bigl)  ^{d-1}; \label{eq:BGoft}
\end{equation}
in this expression $\Omega_{d-1}$ is the area of a unit $d$-1-sphere, and we define 
\begin{equation}
\operatorname{sn}(k,t)  \equiv 
\left\{ \begin{array}{ccccl} 	
\frac{\sin(\sqrt{k}t)}{\sqrt{k}} & \textrm{for} & k>0& \& & 0 \leq t \leq \frac{\pi}{\sqrt{k}}\\
0 & \textrm{for} & k>0& \& & \ \, \ \ \  \  t \geq  \frac{\pi}{\sqrt{k}}\\
t & \textrm{for} & k=0& \& & \ \ \ \ \ \, t \geq 0\\
 \frac{\sinh(\sqrt{-k}t)}{\sqrt{-k}} & \textrm{for} & k<0& \& & \ \ \ \ \ \, t \geq 0  \ . 
 \end{array} \right.   \label{eq:definitionofsnfunction}
\end{equation}
The BG theorem is that this upperbounds the volume of a geodesic ball in the original space,  
\begin{equation}
\textrm{BG bound}: \ \ \textrm{volume}(t) \ \leq \  \textrm{BG}(t) \  .  \label{eq:astatementofBGbound}
\end{equation}
(This was proved by Bishop \cite{BishopGromov,BishopGromov2}, Gromov contributed monotonicity inequalities.) 

The Bishop-Gromov bound applies at all times, but the inequality is often very far from tight. For example, consider the product of $n$ hyperbolic 2-planes $(\mathbb{H}^{2})^{n}$. This is an Einstein space, which means the Ricci tensor is the same in all directions, $\mathcal{R}_{\mu \nu} = -  g_{\mu \nu}$. It is straightforward to calculate 
the late-time volume growth, which is dominated by the diagonal directions that point approximately equally down every $\mathbb{H}^{2}$, and is given by 
\begin{eqnarray}
\textrm{volume}_{(\mathbb{H}^{2})^{n}}(t) &\sim& \exp \left[ \sqrt{ n   } t  \right] \ . \label{eq:actualforproductofhyperboles}
\end{eqnarray}
To calculate the late-time BG bound for this space, we should consider the $2n$-dimensional hyperbolic space that is stretched to have curvature length $\sqrt{2n-1}$ and therefore Ricci curvature $\mathcal{R}_{\mu \nu} = -  g_{\mu \nu}$. The late-time volume growth for this space is 
\begin{eqnarray}
\textrm{BG}(t) = \textrm{volume}_{\mathbb{H}^{2n}}(t)&\sim& \exp \left[ \sqrt{ 2n-1  } t  \right] \ . \label{eq:BGBGBG}
\end{eqnarray}
There is a large difference in exponents between the actual late-time volume Eq.~\ref{eq:actualforproductofhyperboles} and the BG bound Eq.~\ref{eq:BGBGBG}: the BG bound is nowhere near tight. 

 In this paper we will attempt to diagnose why the BG bound so flagrantly fails to be tight, and we will find that in many cases, including this one, the discrepancy is caused by \emph{shear}. We will discuss the cause and effects of shear in detail in Sec.~\ref{subsec:shearbad}. We will then discuss how we can better incorporate the effect of shear in order to tighten the BG bound. Our investigations will culminate in a new bound, Eq.~\ref{eq:bigmainresulthere}. This new bound is more specialized than the BG bound insofar as it only applies at late times (not all times like the BG bound), it only applies to homogeneous metrics, and (for technical reasons) it only applies to spaces with exclusively non-positive sectional curvatures. For all such spaces, when the sectional curvatures of sections that have one leg down the $X_\textrm{min}^\mu$-direction are unequal, which is to say whenever the traceless transverse tensor
 \begin{equation}
 W_{ab} (X)\equiv   \mathcal{R}_{a\mu \nu b}X^{\mu} X^{\nu} + \frac{1}{d-1}  \mathcal{R}_{\mu \nu}X^\mu X^\nu (g_{ab} - X_a X_b) \ \label{eq:definitionofW}
\end{equation}
is nonzero in the $X^\mu = X_\textrm{min}^\mu$ direction, Eq.~\ref{eq:bigmainresulthere} gives an upperbound on the late-time volume growth that is parametrically tighter than the Bishop-Gromov bound.
The improvement is in the exponent of late-time volume growth. \\

\noindent The bound proved in this paper carries some technical restrictions: it only applies at late time, it only applies to homogeneous metrics, and it only applies to spaces with non-positive sectional curvature. However, these are merely technical conveniences, adopted in the interests of most straightforwardly establishing a proof-of-principle that a better-than-Bishop-Gromov bound can be derived by incorporating the effects of shear. In Sec.~\ref{subsec:enhancements} we will describe  how each of these restrictions can be relaxed. \\

\noindent The Bishop-Gromov bound has been used to prove quantum computational complexity lowerbounds \cite{Brown:2021euk} in the context of the geometrical approach to complexity pioneered by Nielsen \cite{Nielsen1,Nielsen2}. The initial motivation for this paper and its companion \cite{BG1} was to strengthen these bounds by enhancing the Bishop-Gromov theorem, but the enhancements we find are too fine-grained to meaningfully improve the quantum complexity lowerbounds. However, our estimates substantially improve volume-growth estimates for classical examples, including the left-invariant metrics on nil-potent Lie groups studied by Heintze \cite{Heintze} and Milnor \cite{milnor}.

\subsection{Deriving the BG bound from Raychaudhuri} \label{subsec:derivingBG}
To see what the Bishop-Gromov bound is missing---and what we can add back---it'll first be helpful to rederive the BG bound starting with the Raychaudhuri equation. The Raychaudhuri equation \cite{Raychaudhuri:1953yv} describes the evolution of a congruence of geodesics. The derivatives of the tangent vector field  $X^{\mu}$  can be decomposed according to their trace, traceless symmetric, and antisymmetric components as 
\begin{eqnarray}
\textrm{expansion } \ \ \  \ \  \theta & \equiv & \nabla_\mu X^{\mu}  \label{eq:defoftheta} \\
\textrm{shear } \  \ \sigma_{\mu \nu} & \equiv &    \frac{1}{2} ( \nabla_\mu X_\nu +\nabla_\nu X_\mu   - \frac{2}{d-1} (\nabla_a X^a) (g_{\mu \nu} - X_\mu X_{\nu}) ) \label{eq:definitionofshear} \\ 
\textrm{vorticity } \ \ \omega_{\mu \nu} & \equiv &  \frac{1}{2} \left( \nabla_\mu X_\nu - \nabla_\nu X_\mu  \right) \ . \label{eq:defvoriticity}
\end{eqnarray} 
With this decomposition, the $d$+0-dimensional Raychaudhuri equation, derived in Appendix~\ref{appendix:derivationofRay}, is 
\begin{equation}
\theta' = - \frac{1}{d-1} \theta^2 -  \sigma^2 + \omega^2 
- \mathcal{R}_{\mu \nu}X^{\mu} X^{\nu} \, , \label{eq:Raychaudhuri}
\end{equation}
where the derivative is taken along a geodesic $\theta' \equiv X^\mu \nabla_\mu \theta$ and  $\sigma^2 \equiv \sigma_{\mu \nu} \sigma^{\mu \nu}$ and $\omega^2 \equiv \omega_{\mu \nu} \omega^{\mu \nu}$. 
The $d+1$-dimensional version of this equation was used by Penrose and Hawking \cite{Penrose:1964wq,Hawking:1969sw} to prove their black hole singularity theorems, first by relating $\mathcal{R}_{\mu \nu}$ to the stress-tensor $T_{\mu \nu}$ via Einstein's equations, and then using the Raychaudhuri equation to upperbound the expansion of lightsheets. We will instead use Eq.~\ref{eq:Raychaudhuri} to upperbound the volume of geodesic balls.  
For geodesic balls, the vorticity is zero, since there is no angular momentum. Putting $\omega^2 = 0$ and using $- \mathcal{R}_{\mu \nu}X^{\mu} X^{\nu}   \leq  
- \mathcal{R}_{\mu \nu}X_\textrm{min}^{\mu} X_\textrm{min}^{\nu}$ gives 
 \begin{equation}
\theta' +  \frac{1}{d-1} \theta^2   + \sigma^2 \ \leq \ 
- \mathcal{R}_{\mu \nu}X_\textrm{min}^{\mu} X_\textrm{min}^{\nu} \, . \label{eq:preBGbound}
\end{equation}
Since  $\sigma^2 \geq 0$, dropping the shear term results in the inequality
 \begin{equation}
\theta' +  \frac{1}{d-1} \theta^2  \ \leq \ 
- \mathcal{R}_{\mu \nu}X_\textrm{min}^{\mu} X_\textrm{min}^{\nu} \, . \label{eq:preBGbound2}
\end{equation}
(This is the step whose inefficiency we will spend much of the rest of the paper trying to ameliorate.)
Next we note that 
\begin{equation}
\textrm{lemma}: \ \theta_1(0) \geq \theta_2(0) \ \  \& \  \  \forall t, \theta'_1(t) +  \frac{1}{d-1} \theta_1^2 \geq  \theta'_2(t) +  \frac{1}{d-1} \theta_2^2(t)  \ \ \rightarrow \ \  \forall t, \theta_1(t) \geq \theta_2(t) \ , \label{eq:introlemma}
\end{equation} 
or in words you maximize $\theta(t)$ by maximizing the right-hand-side of Eq.~\ref{eq:preBGbound2}. 
(For more discussion of this and related lemmas, see Sec.~2.2 of \cite{BG1} using the translation $\theta \equiv \frac{d}{dt} \log j$.)  Finally, integrating Eq.~\ref{eq:preBGbound2} along every geodesic departing the origin and appealing to Eq.~\ref{eq:introlemma}  yields the BG bound, Eq.~\ref{eq:astatementofBGbound}.

\subsection{Why isn't the BG bound tight?}

Except for maximally symmetric spaces, the Bishop-Gromov bound is never tight. Let's discuss why not. 
\begin{enumerate}
\item Anisotropic $\mathcal{R}_{\mu \nu}$. 

The Bishop-Gromov bound treats all directions leaving the origin as though they are leaving in the worst-case most-expansive direction, $X^\mu_\textrm{min}$. For spaces that are not Einstein spaces (i.e.~for spaces that have anisotropic Ricci tensors so that the value of $\mathcal{R}_{\mu \nu}X^\mu X^\nu$ depends on $X^\mu$), this introduces inefficiency.   This inefficiency is most pronounced at early times, and was partially remedied by the improvement we introduced to the BG bound in a previous paper \cite{BG1}.

\item Cut loci. 

Regions that have been swept out by geodesics more than once still contribute only once to the volume. But the derivation of the BG bound described in Sec.~\ref{subsec:derivingBG} doesn't notice that the region has already been included, and so overcounts the volume. For example, in a compact hyperbolic space the BG bound continues to grow just as fast as it would in a hyperbolic space, even though the volume is bounded. 
\item Turning. 

When $\mathcal{R}_{\mu \nu} X^{\mu} X^{\nu}$ is not conserved along geodesics, i.e.~when $X^\rho X^\mu X^\nu \nabla_\rho \mathcal{R}_{\mu \nu} \neq 0$, then the direction of greatest acceleration is shared between different geodesics. We call this `turning' because the Ricci quadratic form is turning relative to the direction of propagation of the geodesic.  For reasons discussed in \cite{BG1}, sharing is inefficient, and makes the volume grow more slowly. In short this is because geodesic flow conserves phase-space volume, so for every geodesic that turns from a slow-growing  direction to a fast-growing direction there must be another than turns from fast to slow, and we showed that the amount of volume gained by the previously slow-growing geodesic speeding up does not fully compensate for the amount of volume lost by the previously fast-growing geodesic slowing down.  This result relied on a property of correlated Jacobi equations that we call the `additive shuffling lemma' and review in Sec.~\ref{subsec:multiadditive}. Turning thus makes the BG bound not tight. (Since the definition of a `locally symmetric space' is that $\nabla_a \mathcal{R}_{\mu \nu \rho \sigma } = 0$, there is no turning for symmetric spaces.)

\item Shear. 

In the derivation of the BG bound, we dropped the shear term $\sigma^2$. As we saw in Eq.~\ref{eq:preBGbound},  the effect of shear is to make the rate of growth of volume  slower than the limit set by the BG bound. 
An example of a space that gives rise to shear is $(\mathbb{H}^{2})^{n}$, considered in Eq.~\ref{eq:actualforproductofhyperboles}. Since this space is an Einstein metric (and indeed a symmetric space) with neither cut loci nor turning, the entire reason this example fails saturate the BG bound is shear.  As we will see in the next subsection, the reason this space gives rise to shear is that its sectional curvatures are unequal. Specifically, the sectional curvature of sections that lie \emph{within} a single $\mathbb{H}^{2}$ are negative, whereas the sectional curvature of sections that lie \emph{across} two $\mathbb{H}^{2}$s are zero.
\end{enumerate}

\subsection{Why shear impedes volume growth} \label{subsec:shearbad}
Let's give a heuristic description of why shear impedes volume growth. We'll start with a geometrical picture of what shear is. Consider zooming in on a patch of the expanding geodesic sphere. Each point on the sphere has a normal vector $\vec{X}$, pointing out like a porcupine. Consider two neighboring geodesics, separated by a vector $\vec{S}$. Neighboring geodesics may be moving farther apart, or closer together, as the sphere expands, and their relative velocity is determined by the gradient of the normal vector field,  $\dot{S}^\mu = S^\nu \nabla_\nu X^\mu$. `Shear' is the traceless symmetric part of $ \nabla_\nu X_\mu$, which means it is nonzero when the relative velocity of neighboring geodesics depends on the direction in which they are separated. For example, shear is when two geodesics separated by a small step in longitude are moving apart with a different velocity than two geodesics separated by a small step in latitude, so a small circle drawn on the geodesic sphere gets deformed into an ellipse as the sphere expands. 

There are $d-1$ directions that lie within the surface of the sphere, and so there are $d-1$ different sections that have have one leg down the normal direction $X^\mu$ and one leg within the sphere. Eq.~\ref{eq:definitionofshear} tells us that the shear $\sigma^2 \equiv \sigma_{\mu \nu}\sigma^{\mu \nu}$ is  the \emph{variance} of the rates at which these sections are expanding. (By contrast the `expansion' $\theta$ is determined by the \emph{average} rate at which these sections are expanding; it tells us about the rate at which a small circle drawn on the geodesic sphere changes area.)  Thus if all the sections are expanding at the same rate, the shear is zero, whereas if they are expanding at wildly different rates (perhaps some can be expanding rapidly while others are contracting), the shear will be large. 

The rate of change of area of a patch is determined solely by the expansion, $\theta$, but the derivative of the expansion---and so the second derivative of the area of the patch---is affected by the shear. This effect is captured by the Raychaudhuri equation, Eq.~\ref{eq:Raychaudhuri}, but to picture what it going on, consider the equation of geodesic deviation.  The equation of geodesic deviation tells us that the rate at which two neighboring geodesics accelerate apart---the rate at which their relative velocity increases---is controlled by the sectional curvature. Specifically, if $S$ is the separation vector between two neighboring geodesics that both point in the $X$ direction, then
\begin{equation}
\ddot{S}^a = \kappa^{a}_{ \ d}(X)  S^d \ \ \textrm{where} \ \ \kappa^{a}_{ \ d}(X)  \equiv \mathcal{R}^{a}_{\ b c d} X^b X^c \   . \label{eq:geodev}
\end{equation} 
Here the derivative is taken along the geodesic, $\ddot{S}^a \equiv X^\mu   \nabla_\mu  X^\nu \nabla_\nu  S^a$, and we see  $\kappa^{a}_{\ d}$ is minus the matrix of sectional curvatures.  Thus shear (different sections having unequal velocities) is sourced by unequal sectional curvatures (different sections having different accelerations). 

In general the sectional curvature of a  pair of neighboring geodesics will not be a constant of motion, even for homogeneous spaces and even when the separation vector points down a principal axis of 
$\kappa^{a}_{\ d}$ and even when the eigenvalues of $\kappa^{a}_{\ d}$ \emph{are} a constant of motion. This is because $S^\mu$ may `precess' around $X^\mu$, so that $S^{a}\kappa_{ab}S^b/S_cS^c$ is not conserved. However, for special cases, including so-called `symmetric spaces' that have $\nabla_\mu \mathcal{R}_{abcd} = 0$, there is no precession and  
the sectional curvature is conserved. Let's first explore the volume growth in such spaces, and then ask what difference precession makes.

\subsubsection{Conserved sectional curvature $\rightarrow$ shear $\rightarrow$ slow volume growth}  \label{subsec:symmetricspaces}
For symmetric spaces, the sectional curvature is conserved. That means that (as well as there being no `turning'), there is no precession. Consequently, if the separation vector between two neighboring geodesics starts off pointing down a principal axis of 
 $\kappa_{ad} \equiv  \mathcal{R}_{abcd} X^b X^c$, it will remaining pointing down that direction as the sphere expands. The $d-1$ eigenvalues of $\kappa_{ad}$, which we will call $\{\kappa_1, \kappa_2, \ldots \kappa_{d-1}\}$, are minus the principal sectional curvatures of the $d-1$ principal sections that contain a geodesic pointing in the $X$-direction. (For illustrative purposes, let's assume all these sectional curvatures are negative, so all the $\kappa_i$s are positive.) A pair of geodesics separated in the $s_i$-direction have their relative acceleration determined by minus the sectional curvature $\kappa_i$, so at late times the separation grow like 
\begin{equation}
s_i \sim \exp [ \sqrt{\kappa_i}  t ] \ . 
\end{equation} 
The area element in the $X$-directions grows like the product of all the orthogonal separations $s_1 s_2 \ldots s_{d-1}$, so the area element can grow no faster than  
\begin{eqnarray}
\textrm{area element} 
&\sim& \exp \left[  ( \sqrt{\kappa_1}  + \sqrt{\kappa_2}  + \ldots \sqrt{\kappa_{d-1}} ) t \right] \label{eq:squarerootofsectionalcurvatures}  \ .
\end{eqnarray}
Crucially, the sum here is \emph{outside} the squareroot. As we will discuss below, this gives a tighter-than-BG bound for anisotropic $\kappa_i$. This stronger bound is tight for the $(\mathbb{H}^{2})^{n}$ example from Eq.~\ref{eq:actualforproductofhyperboles}, and indeed (up to considerations of cut loci) it is tight for all symmetric spaces \cite{symmetricspaces}. 
From the perspective of the Raychaudhuri Equation, the reason we are able to make a  powerful bound on the volume growth when the sectional curvature is conserved is because when the sectional curvature is conserved we can fully account for the effect of the shear. Let's see what changes when we turn on precession.

\subsubsection{Perfect precession $\rightarrow$ no shear $\rightarrow$ faster volume growth} 

When $\nabla_\mu \mathcal{R}_{abcd} \neq 0$, the sectional curvature $S^{a}\kappa_{ab}S^b/S_cS^c$ is generally not conserved. Even if the spectrum $\{\kappa_1, \kappa_2 \ldots \kappa_{d-1}\}$ is preserved by geodesic flow, a pair of geodesics that are initially accelerating apart with acceleration $\kappa_1$ may `precess' and come to accelerate apart with acceleration $\kappa_2$ (since geodesic flow conserves phase space, at the same time a different pair must come to accelerate with acceleration $\kappa_1$). In the last subsection we considered what happens when there is no precession. Now let's consider the opposite limit---precession is so rapid that all pairs of geodesics effectively accelerate apart with the average acceleration (as we will carefully show in Appendix~\ref{appendix:multiplicativeinequalitiesforJacobi}.) In this limit, the separation of any pair of geodesics grows like 
\begin{equation}
s_i \sim \exp \Bigl[ \sqrt{ \frac{\kappa_1 + \kappa_2 + \ldots \kappa_{d-1}}{d-1} }  t \Bigl] \ ,
\end{equation} 
and so the area element (which is the determinant of the $d-1$ separation vectors) grows like 
\begin{equation}
\textrm{area element} 
\sim \exp \Bigl[ \sqrt{(d-1) ({\kappa_1 + \kappa_2 + \ldots \kappa_{d-1}})}  t \Bigl] \ . \label{eq:perfectprecession}
\end{equation}
Unlike in Eq.~\ref{eq:squarerootofsectionalcurvatures}, the sum is now \emph{inside} the squareroot. Since the sum of the sectional curvatures is the trace of the Riemann tensor, and since the trace of the Riemann tensor is the Ricci curvature,  perfect precession recovers the Bishop-Gromov bound Eq.~\ref{eq:BGoft}. 

\subsubsection{Precession makes volume grow faster} 
 By comparing Eq.~\ref{eq:squarerootofsectionalcurvatures} and Eq.~\ref{eq:perfectprecession}, we can see the effect of precession, and therefore of shear. The point is that 
\begin{equation}
\sqrt{(d-1) ({\kappa_1 + \kappa_2 + \ldots \kappa_{d-1}})} \ \geq \ \sqrt{\kappa_1}  + \sqrt{\kappa_2}  + \ldots \sqrt{\kappa_{d-1}} , 
\end{equation}
with equality only when all the $\kappa_i$ are equal. In other words, the effect of precession is to make the volume grow faster. This is because it is more efficient to share the sectional curvatures of greatest acceleration between all the pairs of neighboring geodesics. Hoarding the largest sectional curvature for one pair of geodesics makes them accelerate apart rapidly, but it's not worth it for the overall volume growth because it forces another pair to accelerate apart slowly. The Bishop-Gromov bound thus embodies the `worst case' (i.e.~fastest growth) assumption that precession is perfect and all neighboring geodesics accelerate apart with the same effective sectional curvature---the average sectional curvature---and therefore that the shear is zero. 

\subsubsection{Finite precession $\rightarrow$ bounded shear $\rightarrow$ our new bound}
When precession if perfect, the shear is zero and we get the Bishop-Gromov bound. When precession is zero, we can exactly account for the effect of shear, and we get the tighter bound in Ref.~\cite{symmetricspaces}. This paper is about what happens when you have intermediate amounts of precession, and therefore intermediate amounts of shear. We will develop a bound that is of intermediate strength---stronger than the BG bound, but weaker than the no-precession bound in \cite{symmetricspaces}. (And unlike \cite{symmetricspaces}, our bound applies even for $\nabla_\mu \mathcal{R}_{abcd} \neq 0$.)

Our bound has an advantage over the late-time BG bound whenever the most expansive direction in the BG bound Eq.~\ref{eq:BGoft} has unequal sectional curvatures, which is to say whenever the $W_{ab}$ tensor from Eq.~\ref{eq:definitionofW} is nonzero in the $X_\textrm{min}$-direction. In that case, there will be shear generated in the $X_\textrm{min}$-direction, and by lowerbounding the average shear we will provide a tighter bound on the average volume growth. But as we have just seen, even if the sectional curvatures are unequal, the effect of their inequality can be softened by precession. Part of our task will therefore be upperbounding the rate of precession. This will necessitate not only considering $W_{ab}$ (which must be nonzero to generate shear and give us the opportunity to improve on the BG bound), but also $X^{\mu} \nabla_{\mu} W_{ab}$ (which must be non-infinite to upperbound the rate of precession and therefore lowerbound the rate of shear). 

As we saw, the quantity $\mathcal{R}_{\mu \nu}X^{\mu} X^{\nu}$ is the same for the two spaces $(\mathbb{H}^{2})^{n}$ and $\mathbb{H}^{2n}$, when suitably scaled, and so no bound that only knows about $\mathcal{R}_{\mu \nu}X^{\mu} X^{\nu}$ can possibly capture the fact that the expansion rates of the two spaces are different\footnote{This is a special case of the general principle that if a bound vol$_g(t)$ only knows about geometric data $\{g\}$, it must grow at least as fast as the fastest growing metric with the same data $\{g\}$, from which perspective fast-growing metrics are `barriers' to proving volume upperbounds.}. Our new bound must thus include additional geometrical data beyond the BG quantity, $\mathcal{R}_{\mu \nu}X^{\mu} X^{\nu}$, and the new data it includes is the tensors $W_{ab}$ and $X^\mu \nabla_{\mu} W_{ab}$, as we will now see.  

\subsection{Our new bound} \label{subsec:principalresult}
The principal result of this paper, which we will prove in Sec.~\ref{sec:2}, is that---for $d$-dimensional homogeneous spaces with non-positive sectional curvatures---the time-averaged expansion $\langle \theta \rangle$ is upperbounded by
\begin{equation}
\boxed{ \frac{\langle \theta \rangle^2}{d-1}   \leq    \textrm{max}_{\vec{X} } \Biggl[  R(X)^2 
 - \left( \frac{ \sqrt{4 d R_\textrm{max}^2  {{W}(X)^2}   +  (W_\textrm{max}')^2 } -
   \sqrt{  (W_\textrm{max}')^2  }}{2 d R_\textrm{max}^2  }\right)^2  \Biggl]  } \ . \label{eq:bigmainresulthere} 
\end{equation}
The time-averaged expansion sets the exponent of the late-time volume growth, since, as we will discuss, $\textrm{volume} \leq \exp[ \langle \theta \rangle t ]$.
The terms in this expression as defined by
\begin{eqnarray}
R(X)^2 &\equiv& -\mathcal{R}_{\mu \nu} X^\mu X^\nu\\
W(X)^2 & \equiv & \textrm{Tr}[\hat{W}^2(X)] \equiv W_{ab}(X)W^{ab}(X)  \\
(W'_\textrm{max})^2 & \equiv &  \textrm{max}_{\vec{X}} \textrm{Tr}[(X^c \nabla_c \hat{W}(X))^2] \equiv \textrm{max}_{\vec{X}}  [ (X^c \nabla_c W_{ab}(X))( X^d \nabla_d W^{ab}(X)) ] , 
\end{eqnarray}
where $W_{ab}(X) \equiv \mathcal{R}_{a\mu \nu b}X^{\mu} X^{\nu} + \frac{1}{d-1}  \mathcal{R}_{\mu \nu}X^\mu X^\nu (g_{ab} - X_a X_b)$ introduced in Eq.~\ref{eq:definitionofW} is the traceless part of  ${\kappa}_{ab}$, and $R^2_\textrm{max} \equiv \textrm{max}_{\vec{X}} R(X)^2$, and $X$ is  a unit vector.  
This bound is a function of the geometry only, and depends only on the curvature and its first derivative at a point: to evaluate this bound we do not need to be able to solve the geodesic equation. 

If the direction $X^\mu$ that maximizes $R^2$ gives a nonzero  ${W}_{ab}$, then Eq.~\ref{eq:bigmainresulthere} gives a bound  on $\langle \theta \rangle$ that is strictly tighter than the BG bound, $\langle \theta \rangle^2 < (d-1)R_\textrm{max}^2$. (By contrast for ${W}_{ab}= 0$,  we recover the time-averaged BG bound exactly.) Let's prove Eq.~\ref{eq:bigmainresulthere} now.

\section{Proving our new bound} \label{sec:2}
In Sec.~\ref{subsec:symmetricspaces} we showed that \emph{if the sectional curvature is a constant of motion}, then a bound can be placed on the volume growth that is much more powerful than the Bishop-Gromov bound, because we can fully account for the effect of shear \cite{symmetricspaces}. When the sectional curvature is \emph{not} a constant of motion, the best we could do was the Bishop-Gromov bound. In this section, we'll go beyond the  BG bound. We'll first upperbound the average rate of change of sectional curvature, and then use that to  upperbound the average volume growth. 

To derive the bounds, we will need to track how the expansion $\theta$ and shear ${\sigma}_{ab}$ evolve under geodesic flow. The evolution of $\theta$ is described by the first Raychaudhuri equation, given in Eq.~\ref{eq:Raychaudhuri}. The evolution of ${\sigma}_{ab}$ is described by the second Raychaudhuri equation, 
\begin{equation}
{\sigma}'_{ab}  + \frac{2 \theta \sigma_{ab}}{d-1} + \sigma_{a\mu}\sigma^{\mu}_{\ b} -\frac{\sigma_{\mu \nu} \sigma^{\mu \nu} h_{ab}}{d-1} = W_{ab} \ .  \label{eq:rayforshearmyversion1}
\end{equation}
The matrix $W_{ab}$, which we defined in Eq.~\ref{eq:definitionofW}, is the traceless transverse part of $\kappa_{ab}$, and 
\begin{equation}
h_{ab} \equiv g_{ab} - X_a X_b
\end{equation} is the transverse part of the metric. The first and second Raychaudhuri equations are well-known geometrical identities \cite{Raychaudhuri:1953yv}, but for completeness we rederive them in Appendix~\ref{appendix:derivationofRay}. \\

\noindent For technical reasons, we will restrict ourselves to homogeneous spaces for which every sectional curvature is non-positive. We believe both that these technical restrictions can be relaxed, and  that the bound can be further optimized, as discussed in Sec.~\ref{subsec:enhancements}. However, our analysis will serve as a proof-of-principle that it is possible to improve on the Bishop-Gromov bound even for non-symmetric spaces.

\subsection{Long-time averaged quantities} \label{sec:longtermaverages}
Unlike the original Bishop-Gromov bound, which applies at all times, our new bound only applies at late times. It bounds the exponent of the time-averaged volume growth. Let's discuss why we have to reconcile ourselves to considering only time-averaged quantities. 
\subsubsection{Motivation: why we need to consider time-averaged quantities} 
Recall that the expansion $\theta$ in a particular direction is the rate of change of the logarithm of the area element, so the growth of the area element in that direction is the  exponential of the time-average of the expansion, 
\begin{equation}
\theta \equiv \frac{\partial_t (\textrm{area element})}{ \textrm{area element}} \ \  \rightarrow  \ \ \textrm{area element} \sim \exp[ \langle \theta \rangle t ] \ . \label{eq:averageareaexpansion}
\end{equation}
Now take the $\langle$time-average$\rangle$ of the Raychaudhuri equation, Eq.~\ref{eq:Raychaudhuri}. So long as $\theta$ stays bounded (which will, as we will discuss below, involve ensuring that the geodesic doesn't encounter any caustics and also starting our analysis when the geodesic sphere is of small but nonzero size), the long-term time average allows us to drop the $\langle \theta' \rangle$ term. This leaves  
\begin{equation}
\frac{1}{d-1} \langle \theta^2 \rangle =  \langle -  \mathcal{R}_{\mu \nu}X^{\mu} X^{\nu} \rangle  - \langle  \sigma^2  \rangle. \label{eq:timeaverageRay}
\end{equation}
As we have discussed, the  time-averaged BG bound can then be recovered by using $\langle \sigma^2 \rangle \geq 0$.
But in this paper we will not settle for setting $\langle \sigma^2 \rangle$ to zero. Instead, our goal  will be to lowerbound $\langle \sigma^2 \rangle$ above zero, so that we can then improve on the BG bound. \\

\noindent The reason we have to consider $\langle \sigma^2 \rangle$ and not $\sigma^2$ is that we are {not} going to be able to lowerbound $\sigma^2$ above zero, since it may momentarily be true that $\sigma^2 = 0$ along the geodesic. (For example, at $t=0$ a geodesic ball will always have $\sigma^2 = 0$.) However,  if the shear is zero at one moment, it will not in general be zero at the next. If at some moment the shear is zero $\hat{\sigma} =0 \leftrightarrow \textrm{Tr}[\hat{\sigma}^2] \equiv \sigma^2 \equiv \sigma_{\mu \nu} \sigma^{\mu \nu}= 0$, then Eq.~\ref{eq:rayforshearmyversion1}  reduces to 
\begin{equation}
{\sigma}'_{ab} \biggl|_{\hat{\sigma} = 0} = W_{ab} \ . 
\end{equation}
(The $\hat{\textrm{hat}}$ reminds us that a quantity, e.g.~$\hat{\sigma}$,  is a matrix not a number; we also sometimes instead explicitly write the indices $\sigma_{ab}$.) 
For a typical direction in a typical metric $\hat{W}$ will typically be nonzero. This means that for all except a measure zero set of points along a typical geodesic we'll have $\sigma^2 > 0$, which implies that $\langle \sigma^2 \rangle > 0$ and that the growth rate exponent must be slower than BG.

\subsubsection{Long-term averages from Raychaudhuri}
\noindent We can write Raychaudhuri's second equation, Eq.~\ref{eq:rayforshearmyversion1}, as 
\begin{eqnarray}
\hat{\sigma}' + \hat{f} &=& \hat{W} \ , \label{2lkn34} \\
f_{ab}  & \equiv &  \frac{2 \theta \sigma_{ab}}{d-1} + \sigma_{a\mu}\sigma^{\mu}_{\ b} -  \frac{\sigmaTWO h_{ab}}{d-1} \label{defintionoff1234}, 
\end{eqnarray}
where $\hat{W}$ is defined by Eq.~\ref{eq:definitionofW}. 
Differentiating Eq.~\ref{2lkn34} gives
\begin{equation}
\hat{\sigma}'' + \hat{f}' = \hat{W}' \ . \label{slknerb}
\end{equation}
Here the derivative is taken along a geodesic, so $\hat{f}' \equiv X^c \nabla_c \hat{f}$, etc. Combining Eqs.~\textcolor{red}{\ref{2lkn34}} and \textcolor{blue}{\ref{slknerb}} gives
\begin{eqnarray}
{\color{blue}( \hat{\sigma}'' + \hat{f}') }\hat{\sigma}  + {\color{red}(\hat{\sigma}' + \hat{f})^2} & =&{  \color{blue} \hat{W}'} \hat{\sigma}  + {\color{red} \hat{W}^2} \\
\rightarrow \left( {\color{blue} \hat{\sigma}''}\hat{\sigma} + {\color{red} \hat{\sigma}'\hat{\sigma}' }\right)  + ({\color{blue}\hat{f}'} \hat{\sigma} + {\color{red} \hat{f} \hat{\sigma}'}) + {\color{red}  \hat{\sigma}' \hat{f}} + {\color{red} \hat{f}^2 }&=&{  \color{blue} \hat{W}'} \hat{\sigma}  + {\color{red} \hat{W}^2} \\
\rightarrow  \frac{d}{dt} (  \hat{\sigma}' \hat{\sigma} +  \hat{f} \hat{\sigma}) + \hat{\sigma}' \hat{f} + \hat{f}^2 &=& \hat{W}' \hat{\sigma} + \hat{W}^2 \ . 
\end{eqnarray} 
Now we take the trace and the $\langle$longterm average$\rangle$ along a geodesic. Since all the quantities are bounded after initial transients have died away (as we will discuss), this kills the total derivatives
\begin{equation}
\langle \textrm{Tr}[  \hat{\sigma}' \hat{f} ] \rangle  + \langle \textrm{Tr}[ \hat{f}^2 ]\rangle = \langle \textrm{Tr}[ \hat{W}' \hat{\sigma} ] \rangle+\langle \textrm{Tr}[ \hat{W}^2 ] \rangle \ .  \label{097234h}
\end{equation}
The first term can be evaluated by starting with Eq.~\ref{defintionoff1234}, eliminating ${\sigma}'$ in favor of $\theta'$ by dropping total derivatives, and then substituting in for $\theta'$ using Eq.~\ref{eq:Raychaudhuri}, 
\begin{eqnarray}
\langle \textrm{Tr}[ \hat{\sigma}' \hat{f} ] \rangle &=& \langle \textrm{Tr}[ \hat{\sigma}'  \frac{2 \theta \hat{\sigma}}{d-1}  ] \rangle + \langle \textrm{total derivatives} \rangle  \\
&= & \Bigl\langle \frac{d}{dt} \left(  \frac{\sigmaTWO \theta}{d-1} \right) -  \frac{\sigmaTWO {\theta'}}{d-1} \Bigl\rangle = - \Bigl\langle \frac{\sigmaTWO {\theta'}}{d-1} \Bigl\rangle  \\
&= &  \Bigl\langle \frac{\sigmaTWO (  \frac{1}{d-1} \theta^2 + \sigmaTWO - R^2)}{d-1} \Bigl\rangle \ .  \label{09nwerlkn}
\end{eqnarray}
In the last line we have defined $R^2 \equiv -\mathcal{R}_{\mu \nu} X^\mu X^\nu$. Starting with Eq.~\ref{097234h}, then substituting for $\langle \textrm{Tr}[ \hat{\sigma}' \hat{f}] \rangle$ using Eq.~\ref{09nwerlkn} and $\hat{f}^2$ using Eq.~\ref{defintionoff1234} gives 
\begin{equation}
\boxed{(d-1)\langle \sigmaFOUR \rangle  + 4 \langle \theta \, \sigmaTHREE \rangle + \frac{5 \langle \theta^2 \sigmaTWO \rangle}{d-1}  - \langle  \sigmaTWO R^2 \rangle = (d-1)\langle \textrm{Tr}[ \hat{W}' \hat{\sigma} +  \hat{W}^2 ] \rangle } \ .  \label{eq:boxedversionofaverageofsigmatofour}
\end{equation}
This equation is exact. The equation clearly cannot be satisfied for $\hat{\sigma} = 0$, $\hat{W} \neq 0$, and so it lowerbounds the average shear.

\subsection{Bounds for non-positive sectional curvatures}  \label{subsec22}
\subsubsection{Non-positive sectional curvature implies positive eigenvalues of $\nabla_\mu X^\nu$}  

One obstacle to converting Eq.~\ref{eq:boxedversionofaverageofsigmatofour} to a bound on $\langle \sigmaTWO \rangle$ is that in general $\theta$ can go negative, which complicates the analysis. If a geodesic passes through a caustic we'll have $\theta \rightarrow - \infty$, and then it'll pop out at $\theta = + \infty$ on the other side, creating accounting difficulties. 

 We will discuss an approach to dealing with negative $\theta$ in Sec.~\ref{subsec:enhancements} but, in pursuit of a minimal viable theorem, for the main result of this paper we will prevent $\theta$ from ever going negative by adopting a restriction on the geometry. This restriction is that we will only consider spaces for which the sectional curvatures are all non-positive. When all the sectional curvatures of the space are non-positive, the curvature only ever pushes geodesics apart, so the expansion $\theta \equiv \nabla_\mu X^\mu$ is always positive, and the area element of the  geodesic sphere never stops growing. Indeed, in Appendix~\ref{eq:provingnegativesectionalcurvatureleadstoexponentialgrowth} we'll show the stronger result that 
\begin{equation}
\textrm{lemma: non-positive sectional curvature } \rightarrow \textrm{ every eigenvalue of } \nabla_\mu X^\nu  \textrm{ is  positive.} \label{eq:lemma}
\end{equation}
Non-positive sectional curvatures means that not only is the trace (the sum of the eigenvalues) of the expansion matrix positive, but every single eigenvalue is individually positive. 

\subsubsection[Bounds on $\theta$ when sectional curvatures non-positive]{Bounds on $\theta$ when sectional curvatures are non-positive} \label{sec:iokneb2rnlkEXPANSION}
We just saw that when the sectional curvature is non-positive, the expansion is lowerbounded by zero, $\theta > 0$. Now let's upperbound the late-time expansion. We cannot upperbound the expansion at all times, since at  $t=0$ the expansion is infinite $\theta = \infty$.  But this is a transient effect. The expansion quickly decays away and the Raychaudhuri equation Eq.~\ref{eq:Raychaudhuri} guarantees that once $\theta$ is smaller than $\sqrt{d-1} R_\textrm{max}$ it will stay smaller forever.
After transients have subsided, 
\begin{equation}
\textrm{at late times:} \ \ 0 \ \leq \ \theta \ \leq \ \theta_\textrm{max} \equiv \sqrt{(d-1)} R_\textrm{max}  \ .  \label{2098234809}
\end{equation}
(We have defined a new quantity $\theta_\textrm{max}$ in terms of $ R_\textrm{max}^2 \equiv \textrm{max} |_{\vec{X}}    (- \mathcal{R}_{\mu \nu}X^{\mu} X^{\nu})$.)

\subsubsection[Bounds on $\sigmaTWO$, $\sigmaTHREE$, $\sigmaFOUR$ when sectional curvatures non-positive]{Bounds on $\sigmaTWO$, $\sigmaTHREE$, $\sigmaFOUR$ when sectional curvatures are non-positive}\label{sec:iokneb2rnlkSHEAR}

Let's use the fact that the matrix $\hat{M} \equiv \nabla_\mu X^\nu$ is a positive symmetric matrix to bound the late-time values of the shear.  In terms of $\hat{M}$, the expansion Eq.~\ref{eq:defoftheta} and shear \ref{eq:definitionofshear} can be defined as $\theta \equiv \textrm{Tr}[\hat{M}] $ and $\hat{\sigma} \equiv \hat{M} - \frac{1}{d-1} \textrm{Tr}[\hat{M}] \, \hat{h}$, so that $\sigmaTWO = \textrm{Tr}[\hat{M}^2] - \frac{1}{d-1} \textrm{Tr}[\hat{M}]^2$. (Recall that $\hat{h}$ is defined as the transverse part of the metric, $h_{ab} \equiv g_{ab} - X_a X_b$.) 
%
%
All matrices have $\textrm{Tr}[\hat{M}^2] \geq  \frac{1}{d-1} \textrm{Tr}[\hat{M}]^2$, but matrices with non-negative eigenvalues also have $\textrm{Tr}[\hat{M}^2] \leq \textrm{Tr}[\hat{M}]^2$, which implies that 
\begin{eqnarray}
\sigmaTWO & \leq & \frac{d-2}{d-1} \theta^2 \ . \label{soinsn}
\end{eqnarray}
Combining with Eq.~\ref{2098234809} this gives 
\begin{eqnarray}
\textrm{at late times:} \ \ 0 \  \leq & \sigmaTWO & \leq \ \frac{d-2}{d-1} \theta^2 \ \leq  \  \sigma_\textrm{max}^2  \equiv (d-2) R_\textrm{max}^2 \ .  \label{eq:203948}
\end{eqnarray}
 Let's also upperbound $\sigmaTHREE$ and $\sigmaFOUR$. The following bounds apply in general, without having to assume anything about the positivity of $\hat{M}$. For general $\hat{M}$ we have by the Cauchy-Schwarz inequality that $\textrm{Tr}[\hat{M}^3]^2 \leq (\textrm{Tr}[\hat{M}^2])^3$ and $\textrm{Tr}[\hat{M}^4] \leq (\textrm{Tr}[\hat{M}^2])^2$; from these it immediately follows that $\textrm{Tr}[\hat{\sigma}^3]^2 \leq (\textrm{Tr}[\hat{\sigma}^2])^3$ and $\textrm{Tr}[\hat{\sigma}^4] \leq (\textrm{Tr}[\hat{\sigma}^2])^2$. But notice that the $\hat{M}$ that saturates the Cauchy-Schwarz inequalities has only one nonzero eigenvalue, which is not consistent with the tracelessness of $\hat{\sigma}$, so we can write down stronger inequalities using the fact that $\hat{\sigma}$ is traceless. These inequalities are saturated when one eigenvalue of $\hat{\sigma}$ is $d-2$ and the other $d-2$ eigenvalues are minus one:
\begin{eqnarray}
\left(\sigmaTHREE \right)^2   & \leq & \frac{(d-3)^2}{(d-1)(d-2)} (\sigmaTWO)^3     \label{eq:sigmacubedone}\\
\sigmaFOUR & \leq & \frac{(d-2)^3 +1 }{(d-1)^2(d-2)} (\sigmaTWO)^2 \ . \label{eq:sigmafourthone} 
\end{eqnarray}

\subsubsection[Bounds on $\langle \sigmaTWO \rangle$ when sectional curvatures non-positive]{Bounds on $\langle \sigmaTWO \rangle$ when sectional curvatures are non-positive} \label{subsubsec222}
Equation~\ref{eq:boxedversionofaverageofsigmatofour} gave us a lowerbound on a quantity that included  $\langle \sigma^4 \rangle$. But this does not by itself lead to a lowerbound on $\langle\sigma^2 \rangle$, since in general $\langle \sigma^4 \rangle$ can be unboundedly larger than $\langle \sigma^2 \rangle^2$ if the distribution is very spiky. To limit how spiky the distribution can be, we will use the upperbounds on $\sigma^2$ from Sec.~\ref{subsubsec222}.  Thus by upperbounding $\sigma^2$ we lowerbound $\langle \sigma^2 \rangle$. 

(For a simplified example of this `seesaw' phenomenon, consider a scalar variable $x(t)$. Knowing that $\langle x^4 \rangle = A^2$ gives an \emph{upper}bound $\langle x^2 \rangle<A$, but not a \emph{lower}bound on $\langle x^2 \rangle$. But knowing that $\langle x^4 \rangle = A^2$ \emph{and} that $x^2  < B$ allows us to use the positivity of $\langle x^2 ( B - x^2)\rangle =   B \langle x^2 \rangle - A^2$ to lowerbound $\langle x^2 \rangle >  A^2/B$.)\\

\noindent We will now upperbound the left-hand side of Eq.~\ref{eq:boxedversionofaverageofsigmatofour} by doing three things to make it larger: 
\begin{itemize}
 \item We  use Eqs.~\ref{eq:sigmacubedone} \& \ref{eq:sigmafourthone} to upperbound $\sigmaFOUR$  and $\sigmaTHREE$ in terms of powers of $\sigmaTWO$. 
 \item We  replace some $\theta$s and $\sigmaTWO$s with their maximum values, given by Eqs.~\ref{2098234809} \& \ref{eq:203948}.
 \item We drop $-\langle  \sigmaTWO R^2 \rangle$, since it is non-positive. 
 \end{itemize} 
With these three changes, the left-hand-side of Eq.~\ref{eq:boxedversionofaverageofsigmatofour} is  thus less than 
\begin{equation}
\frac{d^2 - 5d + 7}{d-2}  \langle \sigmaTWO \rangle \sigma_\textrm{max}^2+ 4 \theta_\textrm{max} \frac{d-3}{\sqrt{(d-1)(d-2)}} \langle  \sigmaTWO \rangle  \sqrt{\sigma_\textrm{max}^2 } + \frac{5 \theta^2_\textrm{max}  \langle \sigmaTWO \rangle}{d-1}  \ . \label{1098234} 
\end{equation}
On the other hand, the right-hand side of Eq.~\ref{eq:boxedversionofaverageofsigmatofour} is, by  Cauchy-Schwarz, no less than 
\begin{equation}
- (d-1) \sqrt{ \langle\textrm{Tr}[( \hat{W}')^2] \rangle}  \sqrt{ \langle \sigmaTWO \rangle} + (d-1) \langle \textrm{Tr}[ \hat{W}^2] \rangle \ . 
\end{equation}
Putting these two halves together, and substituting for $\theta_\textrm{max}$ and $\sigma_\textrm{max}$ in terms of $R_\textrm{max}$, gives 
\begin{equation}
d \langle \sigmaTWO \rangle R_\textrm{max}^2 + \sqrt{ \langle \textrm{Tr}[( \hat{W}')^2] \rangle}  \sqrt{ \langle \sigmaTWO \rangle}  \geq \langle \textrm{Tr}[ \hat{W}^2] \rangle \ . \label{eq:quadraticequationforroottracesigma}
\end{equation}
When $\langle \textrm{Tr}[ \hat{W}^2] \rangle  > 0$, this cannot be satisfied by $\langle \sigmaTWO \rangle = 0$, and so forces $\langle \sigmaTWO \rangle$ to be positive. Equation~\ref{eq:quadraticequationforroottracesigma} is a quadratic equation for $\sqrt{\langle {\sigmaTWO} \rangle}$, taking the positive branch gives 
\begin{equation}
\sqrt{\langle \sigmaTWO \rangle} \ \geq \ 
\frac{ \sqrt{4 d R_\textrm{max}^2 \langle \textrm{Tr}[\hat{W}^2] \rangle + \langle \textrm{Tr}[(\hat{W}')^2] \rangle} -
   \sqrt{  \langle \textrm{Tr}[(\hat{W}')^2] \rangle }}{2 d R_\textrm{max}^2} \ . \label{eq:lowerboundonsigmasquared1234}
\end{equation}
\subsubsection[Bounds on $\langle \theta \rangle$ when sectional curvatures non-positive]{Bounds on $\langle \theta \rangle$ when sectional curvatures are non-positive}
Using the time-averaged Raychaudhuri's equation, Eq.~\ref{eq:timeaverageRay}, the fact that $\langle \theta \rangle^2 \leq \langle \theta^2 \rangle$, and  the bound on  $\langle\sigmaTWO\rangle$ from Eq.~\ref{eq:lowerboundonsigmasquared1234} gives an upperbound on the asymptotic growth rate 
\begin{equation}
\frac{\langle \theta \rangle^2}{d-1} \ \leq \    \langle R^2 \rangle -  \left( \frac{ \sqrt{4 d R_\textrm{max}^2 \langle \textrm{Tr}[\hat{W}^2] \rangle + \langle \textrm{Tr}[(\hat{W}')^2] \rangle} -
   \sqrt{  \langle \textrm{Tr}[(\hat{W}')^2] \rangle }}{2 d R_\textrm{max}^2} \right)^2   \label{098234lkn} 
\end{equation}
Recall that $R$, $\hat{W}$, and $\hat{W}'$ are all functions of the direction $\vec{X}(t)$ in which the geodesic is pointing. We wish to turn Eq.~\ref{098234lkn} into a bound that is a function of the geometry only (so that we don't have to calculate $\vec{X}(t)$). One way to do so is to use the fact that the average of a quantity cannot be larger than its largest value (nor smaller than its smallest value) and replace $\langle R^2 \rangle $ by $R^2_\textrm{max} \equiv  \textrm{max}_{\vec{X}} R^2(X) \equiv \textrm{max}_{\vec{X}} (-\mathcal{R}_{\mu \nu} X^\mu X^\nu)$,  replace $\langle \textrm{Tr}[\hat{W}^2] \rangle$ by $W_\textrm{min}^2 \equiv \textrm{min}_{\vec{X}}\textrm{Tr}[\hat{W}(X)^2]$, and replace $\langle \textrm{Tr}[(\hat{W}')^2] \rangle$ by $(W_\textrm{max}')^2  \equiv  \textrm{max}_{\vec{X}}\textrm{Tr}[(\hat{W}')^2]$, 
to give 
\begin{equation}
\frac{\langle \theta \rangle^2}{d-1} \ \leq \     R_\textrm{max}^2  -  \left( \frac{ \sqrt{4d R_\textrm{max}^2  W_\textrm{min}^2  + (W_\textrm{max}')^2} -
   \sqrt{  (W_\textrm{max}')^2  }}{2  d R_\textrm{max}^2   }\right)^2  \ .   \label{098234lknQQ} 
\end{equation}
Since none of the quantities in this expression are a function of $\vec{X}(t)$, this successfully gives a bound on the time-averaged volume growth that is a function of the geometry only. \\

\noindent Let's make a tighter bound with one small improvement. The worst-case assumption embedded in Eq.~\ref{098234lknQQ} is that the worst-case trajectory spends the whole time pointing in the direction $\vec{X}$ that maximizes $R^2(X)$, and also spends the whole time pointing in the direction $\vec{Y}$ that minimizes $\textrm{Tr}[\hat{W}(Y)^2]$. When $\vec{X}$ and $\vec{Y}$ are different directions, this assumption is impossibly pessimistic. Even the worst-case trajectory can only point in one direction at once. Indeed, if we replace $\langle \textrm{Tr}[(\hat{W}')^2] \rangle$ with $(W_\textrm{max}')^2$ on the right-hand-side of Eq.~\ref{098234lkn}, what remains is the sum of a term linear in $\langle R^2 \rangle$ and a term concave in $\langle \textrm{Tr}[\hat{W}^2] \rangle$. This implies that the worst-case trajectory that maximizes the right-hand-side will be independent of time---the maximum cannot be given by a trajectory that spends part of the time pointing in one direction and part of the time pointing in another. This converts the problem of maximizing over trajectories $\vec{X}(t)$ to the problem of maximizing over directions $\vec{X}$. Consequently we can make a tighter bound by insisting that we choose a single direction to jointly maximize the right-hand side. This allows us to derive the principal result of this paper, Eq.~\ref{eq:bigmainresulthere}. 

({You might wonder why we cannot also force the bound to use the same $\vec{X}$ for $\hat{W}'(X)$ that it uses for $R(X)$ and $\hat{W}(X)$ (i.e.~why you cannot make a stronger bound in which every appearance of $(W_\textrm{max}')^2$ in Eq.~\ref{eq:bigmainresulthere} is replaced by $\textrm{Tr}[\hat{W}'(X)^2]$). This is because Eq.~\ref{098234lkn} is not a concave function of $\textrm{Tr}[\hat{W}'(X)^2]$, and so we are not guaranteed that the trajectory $\vec{X}(t)$ that maximizes the right-hand side of Eq.~\ref{098234lkn} is independent of time.
For a concrete example of what goes wrong, consider two directions in a $d=4$-dimensional geometry, one of which has $R^2(X_1) = \textrm{Tr}[\hat{W}(X_1)^2] =1$, $\textrm{Tr}[\hat{W}'(X_1)^2] =0$, and the other has $R^2(X_2) = 1$, $ \textrm{Tr}[\hat{W}(X_2)^2] = 2$, $\textrm{Tr}[\hat{W}'(X_2)^2] =4$. Simple algebra confirms that the right-hand side of  Eq.~\ref{098234lkn} is larger for the `mixed strategy' trajectory that spends half the time pointing in the $X_1$-direction and half  pointing in the $X_2$-direction than it is for either of the `pure strategy' trajectories that spend all the time pointing in just one of the directions.})

\section{Conclusion and Discussion} \label{section:discussion}

In this paper we proved a new bound on the rate of growth of geodesic balls, Eq.~\ref{eq:bigmainresulthere}. Unlike the Bishop-Gromov theorem \cite{BishopGromov}, this new bound only constrains the \emph{late-time average} rate of growth of volume, and only applies to homogeneous spaces with non-positive sectional curvatures.  When it does apply, however, it can multiplicatively improve on the exponent of the late-time volume growth, as we will see explicitly in the examples of Appendix~\ref{appendix:examplemetrics}. Let's now reflect on how we made this improvement, and where we might seek further improvement. 

\subsection{The meaning of $\hat{W}$ and $\hat{W}'$}
Our new bound, Eq.~\ref{eq:bigmainresulthere}, is a function of three geometrical quantities. The first quantity is the same projected Ricci curvature $\mathcal{R}_{\mu \nu}X^\mu X^\nu$ that is already used in the BG bound; the two new quantities are  $\hat{W}$ and $\hat{W}'$. Let's discuss how the appearance of these quantities is exactly what we should expect given the motivating intuitions we explored in Sec.~\ref{subsec:shearbad}. 

The equation of geodesic deviation, Eq.~\ref{eq:geodev}, tells us that the rate at which neighboring geodesics that point in the $\vec{X}$-direction accelerate apart is governed by the matrix of minus the sectional curvatures, $\kappa_{ab}(X)$. Our new bound is better than the BG bound whenever, in the direction $\vec{X}_\textrm{min}$ that minimizes the Ricci curvature, $\kappa_{ab}$ has unequal eigenvalues. Equivalently we can say that our new bound is better whenever the  traceless part of $\hat{\kappa}(X_\textrm{min})$, 
\begin{equation}
\hat{W} = \hat{\kappa} - \frac{1}{d-1} \textrm{Tr}[\hat{\kappa}] \hat{h} \label{checkthissign}
\end{equation}
is nonzero\footnote{Note that despite similarities $\hat{W}$ is \emph{not} the same as the Weyl tensor projected in the $\vec{X}$-direction,
\begin{equation}
C_{ad} (X) \equiv 
\left(\mathcal{R}_{abcd}  + \frac{\mathcal{R}_{bc} g_{ad} -  \mathcal{R}_{bd} g_{ac}  +  \mathcal{R}_{ad} g_{bc} -  \mathcal{R}_{ac} g_{bd} }{d-2}   + \frac{\mathcal{R}  \left( g_{ac} g_{bd} - g_{ad} g_{bc} \right)  }{(d-1)(d-2)} \right) X^b X^c \ . 
\end{equation}
Though both are symmetric $C_{ab} = C_{ba}$ \& $W_{ab} = W_{ba}$, both are traceless $C^a_{ \ a} = W^a_{ \ a} = 0$, and both are transverse $C_{ab}X^b = W_{ab}X^b = 0$, they differ by an additive factor of the traceless transverse part of the Ricci tensor $\tilde{\mathcal{R}}_{ad} \equiv h_{ab}  \mathcal{R}^{bc}  h_{cd} - \frac{1}{d-1} h_{eb}  \mathcal{R}^{bc}  h_{c}^{\ e} h_{ad}$. For example, though the Weyl tensor is identically zero in three dimensions, the tensor $\hat{W}$ need not be, as we will see with an explicit worked example in Appendix~\ref{subsec:squashedH3}.}. Indeed, we see that the correction term in our bound, Eq.~\ref{eq:bigmainresulthere} is driven by $\textrm{Tr}[\hat{W}^2(X)]$, which is the variance of the sectional curvatures amongst those sections that include the direction $X$. From the perspective of the analysis of Sec.~\ref{subsec:shearbad}, this makes perfect sense. In that section, we saw that a major reason that the BG bound is not tight is that it neglects the effect of shear. Shear is when different sections are expanding at unequal rates, and is therefore sourced by variance in sectional curvatures---it is sourced by $\hat{W}$. Our bound improved on the BG bound by better accounting for the effects of shear, and so our improvement gets bigger with bigger shear, which means it gets bigger with bigger $\hat{W}$. 

Our bound, Eq.~\ref{eq:bigmainresulthere}, also depends on the rate of change of $\hat{W}$ along the geodesic, $X^\mu \nabla_\mu W_{ab}$. But the direction of the dependence is the reverse of what it was for $\hat{W}$: the correction term gets bigger with \emph{smaller} $\hat{W}'$. Again, from the perspective of the analysis of Sec.~\ref{subsec:shearbad}, this is what we would expect. In that section, we saw that $\hat{W}$ drives the growth of shear, but that  $\hat{W}'$ reduces shear. This is because $\hat{W}'$ may make neighboring geodesics `precess', so that the pair that was initially separated in the direction of maximum relative acceleration comes to be separated in a direction with less negative sectional curvature. Since the directions of fastest acceleration are shared amongst many different pairs of neighboring geodesics, the relative velocities are averaged out, and the shear is reduced. In summary, more $\hat{W}'$ means more precession, more precession means less shear, less shear means more volume growth, and more volume growth means the BG bound is closer to tight and the correction term from our new bound is smaller.

\subsection{Multiplicative and additive shuffling} \label{subsec:multiadditive}

In a previous paper \cite{BG1} we proved a completely different enhancement to the BG bound. That enhancement was complementary to the one proved in this paper: the  bound of \cite{BG1} found strongest purchase at early times, whereas the bound in this paper applies only at late times; the bound of \cite{BG1} used the fact that not all geodesics can simultaneously be pointing in the worst-case most expansive direction, whereas the bound in this paper does not; the bound of \cite{BG1} neglected shear, whereas the bound in this paper incorporates some of the effects of shear. There is another sense in which these two papers are complementary, which is that they connect to two complementary lemmas about the growth of correlated families of Jacobi equations: our previous paper used a lemma about `additive shuffling', whereas the analysis of this paper is undergirded by a lemma about `multiplicative shuffling'. 

`Shuffling', a concept we introduced in Ref.~\cite{BG1}, is a permutation one can perform on a correlated family of Jacobi equations. Consider a family of solutions to Jacobi equations $j_i(t)$, indexed by $i$. The solutions all start with $j_i(0) = 0 \ \& \ j'_i(0) = 1$ and obey the equations of motion 
\begin{equation}
j_i''(t) = \kappa_i(t) j_i(t) \ . \label{eq:definitionofjacobi} 
\end{equation}
The $i$th $\kappa$-schedule $\kappa_i(t)$ determines the acceleration of the $i$th trajectory $j_i(t)$. The shuffling lemmas ask what happens to the total growth rate of all the Jacobi solutions in the family if we perform a time-dependent permutation (or `shuffling') of which trajectory $j_i(t)$ is following which $\kappa$-schedule $\kappa_i(t)$. 

The additive shuffling lemma, proved in Ref.~\cite{BG1}, tells us that to maximize the \emph{sum} of the Jacobi solutions, $\sum_i j_i(t)$, we should shuffle the coefficients so as to \emph{maximize} the inequality between the Jacobi solutions. This means we should shuffle the $\kappa_i(t)$ in such a way that one solution is always following whatever is at that moment the largest $\kappa_i(t)$, another solution is always following the second largest, and so on:
\begin{equation}
\textrm{additive shuffling lemma \cite{BG1}: \ \ \ to maximize } \sum_i j_i(t), \textrm{ maximize inequality}. 
\end{equation}  
This was relevant for the bound proved in \cite{BG1}, because it showed that for the purposes of upperbounding the volume growth, we could make the worst-case assumption that one geodesic always followed the most-expansive direction with largest $-\mathcal{R}_{\mu \nu}X^\mu X^\nu$, another the second-most-expansive direction, and so on. Because of the additive shuffling lemma, `turning', which is what we call it when  $\mathcal{R}_{\mu \nu}X^\mu X^\nu$ is not a constant of motion, will only make the volume grow slower, and  can therefore be safely ignored for the purposes of establishing an upperbound. 

By contrast, the results of this paper implicitly relied on the multiplicative shuffling lemma. (Implicitly because this lemma is baked in to the Raychaudhuri equation, so we didn't need to make it explicit.) The multiplicative lemma says that to maximize not the sum but the product of the Jacobi solutions, we should \emph{minimize} inequality, 
\begin{equation}
\textrm{multiplicative shuffling lemma: \ \ \ to maximize } \prod_i j_i(t), \textrm{ minimize inequality}. 
\end{equation}  
This minimization can be achieved by rapid shuffling amongst the $\kappa_i(t)$, so that each trajectory effectively follows the same average $\kappa_\textrm{av.}(t) = \sum_i \kappa_i(t)/\sum_i$ and all the $j_i(t)$ are the same. The multiplicative shuffling lemma is proved in Appendix~\ref{appendix:multiplicativeinequalitiesforJacobi}.  In its connection to the bound proved in this paper, the $j_i(t)$s are the separation of neighboring geodesics, and (in an orthonormal basis) the product of the $j_i(t)$s gives the area element. What this lemma tells us, therefore, is that precession---which averages the sectional curvatures and reduces the inequality in the expansion rates and therefore reduces the shear---makes the area element and therefore the volume grow faster.

We thus see that the shuffling strategy for maximizing the arithmetic mean of the $j_i(t)$s is diametrically opposite to the strategy for maximizing the geometric mean.

\subsection{How to further improve the bound}  \label{subsec:enhancements}
The bound we have established, Eq.~\ref{eq:bigmainresulthere}, is meant to be a proof-of-principle that by incorporating shear we can improve on the  BG bound for the late-time volume growth. In pursuit of the easiest possible such bound, we made a number of simplifying assumptions that limited the power of our theorem. Let's discuss ways to make the theorem better. 

One class of improvements is simply to make our bound tighter by reducing the inefficiency of the steps we took in Sec.~\ref{subsec22}. For example, recall that to reach the upperbound in Eq.~\ref{1098234}, we  deleted the term  $-R^2 \langle  \sigmaTWO \rangle$. If the space has a nonzero  $R_\textrm{min}^2 \equiv \textrm{min}_{\vec{X}} (- \mathcal{R}_{\mu \nu}X^\mu X^\nu)$, we could tighten the bound by instead replacing $-R^2 \langle  \sigmaTWO \rangle$ with $-R_\textrm{min}^2 \langle  \sigmaTWO \rangle$. This would lead to 
\begin{equation}
 \frac{\langle \theta \rangle^2}{d-1}   \leq    \textrm{max}_{\vec{X} } \Biggl[  R^2(X) 
 - \left( \frac{ \sqrt{4 (d R_\textrm{max}^2 - \frac{1}{d-1} R_\textrm{min}^2)  {\textrm{Tr}[\hat{W}(X)^2}]   +  (W_\textrm{max}')^2 } -
   \sqrt{  (W_\textrm{max}')^2  }}{2  (d R_\textrm{max}^2 - \frac{1}{d-1} R_\textrm{min}^2)  }\right)^2  \Biggl]  \ . \label{eq:improvementwithRmin}
   \end{equation}
If the space is an Einstein metric ($R(X) = R_\textrm{min} = R_\textrm{max}$) this simplifies further. 
Similarly, there were a number of other inefficient steps in Sec.~\ref{sec:2} that are ripe for improvement. 

Another class of improvements would be to widen the scope of applicability of our bound. For example, in this paper we restricted our considerations to spaces with exclusively non-positive sectional curvatures. The motivation for this restriction was to stop $\theta$ going negative in order to be able to turn an upperbound on $\theta$ into an upperbound on $\theta^2$, and a similar argument for $\sigma^2$. However, one might seek to remove this restriction. For example, one could attempt to construct an argument that if $\theta$ goes too negative, we can thenceforth neglect that geodesic since it is doomed to end at a caustic and no longer contribute to the volume growth. Recall the Raychaudhuri equation, Eq.~\ref{eq:Raychaudhuri}, 
 \begin{equation}
\theta' = -  \frac{1}{d-1} \theta^2 - \sigmaTWO + R^2 .
\end{equation}
If ever $\theta < - \sqrt{(d-1)} R_\textrm{max}$ this equation guarantees that the negative curvature is not strong enough to turn around the inbound geodesics and stop them reaching $\theta = - \infty$. (This is basically the Penrose-Hawking singularity theorem in the presence of a positive cosmological constant.) Intriguingly, the condition on $\theta^2$ that a geodesic not be contracting so fast that it is doomed to hit a singularity, $\theta^2 < (d-1) R_\textrm{max}^2$, is the same as the bound we derived for positive $\theta$ in Eq.~\ref{2098234809}. A similar argument could upperbound $\sigma^2$ for non-doomed geodesics in terms of $\hat{\kappa}$, allowing one to construct an extended version of our bound that also applies to spaces with positive sectional curvature. 

Another class of improvements would be to widen the definition of the `volume' that we are bounding. According to the definition we adopted at the very start of the paper, the geodesic ball is the union of all points that have been touched at least once. This means that regions that have been swept out twice or more still only contribute to the volume once. But nothing in the proof of Eq.~\ref{eq:bigmainresulthere} required this restriction: the argument was a `local' argument about the surface of the geodesic sphere. Eq.~\ref{eq:bigmainresulthere} thus upperbounds the total volume where we count each point by the number of times it has been swept out, in such a way that it continues to grow even for e.g.~compact hyperbolic spaces.  (If we adopt this more expansive definition of volume we also have a \emph{lower}bound on the volume growth, since $ \theta \geq R_\textrm{min}$.)

In this paper we have considered only the long-term time-averaged growth rate. The motivation for this restriction was so that we could drop total derivative terms in Sec.~\ref{sec:longtermaverages}: a fixed contribution from the total derivative amortized over an infinite timescale gives an average contribution of zero. But this was unduly conservative. We know the maximum possible size of the total derivative terms we dropped, so we can easily bound the total error in our estimates for $\langle \theta \rangle$, and see that it rapidly goes to zero for even intermediate $t$. Bounding the corrections at finite $t$ will  allow us to apply our bound even to spaces of finite total volume such as the deformed unitary groups discussed in \cite{Brown:2021euk}. 

Finally, it would be interesting to widen the applicability of our bound so that it also applies to \emph{inhomogeneous} spaces. For inhomogeneous spaces, the original BG bound requires us to generalize Eq.~\ref{eq:BGoft} so as to take not only the worst-case direction, but the worst-case direction around the worst-case point. It is clear that a version of our bound Eq.~\ref{098234lknQQ} would also apply to inhomogeneous spaces if  each of the quantities in that equation---$R_\textrm{max}$, $W^2_\textrm{min}$ and $(W'_\textrm{max})^2$---is evaluated in the worst-case direction for the worst-case point, taking possibly different directions and different points for each of those three quantities. Thus Eq.~\ref{098234lknQQ} has an immediate generalization to `regular' inhomogeneous spaces.

\section*{Acknowledgements} 
MHF thanks the Aspen Center for Physics for hospitality.

\appendix

\section{The first and second Raychaudhuri equations}  \label{appendix:derivationofRay}
To be self-contained, let's derive the first and second Raychaudhuri equations \cite{Raychaudhuri:1953yv}. 
Recalling that $X^{\mu}$ is the vector field composed of the unit tangents to the geodesics at each point, we will need the geodesic equation $X^{\mu} \nabla_{\mu} X^{\nu} = 0$ and the definition of the Riemann tensor as the commutator of the covariant derivatives $\left( \nabla_\mu \nabla_\nu - \nabla_\nu \nabla_\mu \right) X^\rho = \mathcal{R}_{\ \sigma \mu \nu}^{\rho} X^\sigma \ .$ 
These definitions imply 
\begin{eqnarray}
\frac{d}{d \tau} \nabla_b X^c  = X^a \nabla_a \nabla_b X^c &=& X^a (\nabla_a \nabla_b - \nabla_b \nabla_a) X^c + X^a  \nabla_b \nabla_a  X^c  \\ 
  &=& X^a (\nabla_a \nabla_b - \nabla_b \nabla_a) X^c +   \nabla_b (X^a \nabla_a  X^c)  -   (\nabla_b X^a) (\nabla_a  X^c ) \nonumber \\ 
  & = & \mathcal{R}^c_{\ d a b} X^a X^d -   (\nabla_b X^a) (\nabla_a  X^c ) \ .  \label{eq:prepreRay}
\end{eqnarray}
Alternatively, writing $M_{a}^{\ b} \equiv \nabla_a X^b$ we can write this as 
\begin{equation}
\frac{d \hat{M}}{dt} = \hat{\kappa} - \hat{M}^2 \ , \label{eq:equationofmotionforM}
\end{equation}
where again $\hat{\kappa}$ is minus the matrix of sectional curvatures $\kappa^{c}_{\ b} \equiv  \mathcal{R}^c_{\ d a b} X^a X^d $. 
Raychaudhuri's first equation, Eq.~\ref{eq:Raychaudhuri}, follows by taking the trace of Eq.~\ref{eq:equationofmotionforM}; Raychaudhuri's second equation, Eq.~\ref{eq:rayforshearmyversion1}, follows by taking the traceless symmetric part. 

It is standard to write these equations by breaking down $\hat{M}$ into its irreducible transverse representations. The matrix $\hat{M}$ is `transverse' since $X^{\mu} \nabla_{\mu} X_{\nu} = 0 = X^{\nu} \nabla_{\mu} X_{\nu}$, and we choose the decomposition so that each term individually retains this property and thus lives in the codimension-one subspace orthogonal to $X^{\mu}$; this leads to the definitions Eqs.~\ref{eq:defoftheta}, \ref{eq:definitionofshear} and \ref{eq:defvoriticity} and the form shown in Eqs.~\ref{eq:Raychaudhuri} and \ref{eq:rayforshearmyversion1}. \\

\noindent Though we will not emphasize this perspective in this paper, another away to think about the meaning of $\hat{M}$ is as follows. We can think of it as being \begin{equation}
\hat{M} =  \dot{\hat{T}} \hat{T}^{-1}\ , \label{eq:MintermsofT}
\end{equation}
where $\hat{T}$ is the   transformation matrix that maps the separation at the initial time to a separation at a later time,
\begin{equation}
S^a(t) = T^{a}_{ \  \,b}(t,t_0) S^b(t_0)  \ . 
\end{equation}
The equation of motion for $\hat{T}$ is given by the equation of geodesic deviation, $\ddot{\hat{T}} = \hat{\kappa} \hat{T}$, and then via  Eq.~\ref{eq:MintermsofT} this gives the equation of motion for $\hat{M}$, Eq.~\ref{eq:equationofmotionforM}. This perspective makes it transparent that the expansion is the logarithmic derivative of the area form on the geodesic sphere, $\theta = \textrm{Tr} \hat{M} = \frac{d}{dt} \log [ \det \hat{T}]$.

\section{Negative sectional curvature $\rightarrow$ positive expansion} \label{eq:provingnegativesectionalcurvatureleadstoexponentialgrowth}

The expansion matrix $\hat{M}$ obeys Eq.~\ref{eq:equationofmotionforM}.  Let's prove that if all the eigenvalues of $\hat{\kappa}$ are non-negative, and if at $t=0$ all the eigenvalues of $\hat{M}$ are non-negative, then all the eigenvalues of $\hat{M}$ will remain non-negative forever. This will establish our lemma, Eq.~\ref{eq:lemma}.

 Proving this lemma would be trivial were $\hat{M}$  a one-by-one matrix (a number): it would follow immediately from $\dot{M} + M^2 = \kappa \geq 0$ that $-M^2$ vanishes too fast as $M \rightarrow 0$ to drive $M$ negative. However, for matrix $\hat{M}$ we are going to have to work harder since even when an eigenvalue is zero, the matrix $\hat{M}^2$ is typically nonzero.

Consider the determinant of $\hat{M}$. We will have proved our lemma if we can prove that the determinant cannot cross from positive to negative when $\hat{M}$ and $\hat{\kappa}$ are real symmetric matrices with non-negative eigenvalues. The derivative of the determinant of any matrix is 
\begin{equation}
\frac{d}{dt} \textrm{det}[\hat{M}] = \textrm{Tr}[\textrm{adj}[\hat{M}] . \dot{\hat{M}} ] . 
\end{equation}
Recalling that $\textrm{adj}[\hat{M}] \equiv  \textrm{det}[\hat{M}] \hat{M}^{-1} $ for an invertible matrix, and then applying the equation of motion Eq.~\ref{eq:equationofmotionforM}, gives  
\begin{equation}
\frac{d}{dt} \textrm{det}[\hat{M}] = - \textrm{det}[\hat{M}]\textrm{Tr}[\hat{M} ] + \textrm{Tr}[\textrm{adj}[\hat{M}] . \hat{\kappa} ] \ . \label{lkn234}
\end{equation}
The second term on the right-hand side is positive, since  (i) $\textrm{adj}[\hat{M}]$ is a positive symmetric matrix, since if $\hat{M}$ is positive symmetric then $\textrm{adj}[\hat{M}]$ is positive symmetric; (ii) $\hat{\kappa}$ is positive symmetric by assumption; and (iii) the product of two positive symmetric matrices is a positive (though not necessarily symmetric) matrix. On the other hand, while the first term on the right-hand side is negative, it goes to zero so fast as $\textrm{det}\hat{M} \rightarrow 0$ that it can never drive $\textrm{det}\hat{M}$ negative. One way to see this is to rewrite Eq.~\ref{lkn234} as 
\begin{equation}
\frac{d}{dt} \left( e^{\int\textrm{Tr}[\hat M] dt } \textrm{det}\hat{M}\right)  = e^{\int \textrm{Tr}[\hat M] dt} \textrm{Tr}[\textrm{adj}[\hat{M}] .  \hat{\kappa}  ] \ . 
\end{equation}
Thus for $\hat{\kappa} \geq 0$, we have that $\textrm{det}\hat{M}$ can at worst  dilute away like $e^{-\int \textrm{Tr}[\hat M] dt} = e^{-\int \theta dt}$, and will therefore never hit zero.  This proves that if all the eigenvalues of $\hat{M}$ start positive, and if all the eigenvalues of $\hat{\kappa}$ are forever non-negative, then all the eigenvalues of $\hat{M}$ remain positive. This establishes our lemma, Eq.~\ref{eq:lemma}.

\section{Multiplicative inequalities for Jacobi equations} \label{appendix:multiplicativeinequalitiesforJacobi}

In this appendix, we will consider properties of solutions $j_i(t)$  that obey Jacobi equations Eq.~\ref{eq:definitionofjacobi} and have initial conditions $j_i(0) = 0$ and $j_i'(0) = 1$. We will imagine we are granted the power to permute (or `shuffle') at each moment which trajectory $j_i(t)$ is following which $\kappa$-schedule, and will ask how we should use our power to maximize the product of the $j_i(t)$s.

For the main result of this paper we considered non-positive curvature, which implied non-negative $\kappa$ and therefore positive $j(t)$s.   But we wish to prove a more general lemma that permits negative $\kappa_i(t)$ so that it can also be applied when the curvature is positive, following the path laid out in Sec.~\ref{subsec:enhancements}. If we are considering negative $\kappa$, we may encounter a situation in which the solution to Eq.~\ref{eq:definitionofjacobi} hits zero and potentially goes negative. In that situation we will, for the purposes of our lemma, have an additional \emph{ad hoc} rule:  if the solution to Eq.~\ref{eq:definitionofjacobi} hits zero, then $j_i(t)$ sticks at zero forever thereafter (see \cite{BG1} for more discussion).   
\subsection{Averaging a pair of $\kappa$s increases the product $j_1 j_2$}

Given two trajectories $j''_1(t) = \kappa_1(t) j_1(t)$ and $j''_2(t) = \kappa_2(t) j_2(t)$ we can define a third trajectory $j_\textrm{av.}$ by $j''_\textrm{av.} (t) = \kappa_\textrm{av.} (t) {j}_\textrm{av.} (t)$ with
\begin{equation}
\kappa_\textrm{av.}(t) \equiv \frac{\kappa_1(t) + \kappa_2(t) }{2} \ . 
\end{equation}
Then let's show that, without restriction on the sign of $\kappa$, that for all times 
\begin{equation}
\textrm{multiplicative lemma}:  \  \ \ 2 \frac{j'_\textrm{av.}}{{j}_\textrm{av.}}  \geq  \frac{j'_\textrm{1}}{{j}_\textrm{1}}  +  \frac{j'_\textrm{2}}{{j}_\textrm{2}} \ \  \ \& \  \ \ j_\textrm{av.}^2 \geq j_1 j_2 \ . \label{eq:multiplicativelemma}
\end{equation}
First, notice that these equations are satisfied at $t=0$. (Even though $\frac{j'}{{j}}$ is infinite at $t=0$, the quantity $2 \frac{j'_\textrm{av.}}{{j_\textrm{av.}}} - \frac{j'_1}{{j_1}} - \frac{j'_2}{{j_2}}$ is finite and equal to zero; Taylor expanding gives  $2 \frac{j'_\textrm{av.}}{{j_\textrm{av.}}} - \frac{j'_1}{{j_1}} - \frac{j'_2}{{j_2}} = \frac{1}{45} \left( \kappa_1(0)^2 + \kappa_2(0)^2 -2\kappa_{\textrm{av.}}^2 \right) t^3 + \ldots \geq 0$.) Now let's show that they remain satisfied for all $t$.

\begin{itemize}
\item Argument that it remains true that $2 \frac{j'_\textrm{av.}}{{j}_\textrm{av.}}  \geq  \frac{j'_\textrm{1}}{{j}_\textrm{1}}  +  \frac{j'_\textrm{2}}{{j}_\textrm{2}}$: 

if they were to cross, there would need to be a moment when $2 \frac{j'_\textrm{av.}}{{j}_\textrm{av.}}  =  \frac{j'_\textrm{1}}{{j}_\textrm{1}}  +  \frac{j'_\textrm{2}}{{j}_\textrm{2}}$ and $\frac{d}{dt}  \left( 2  \frac{j'_\textrm{av.}}{{j}_\textrm{av.}} \right)   <\frac{d}{dt}  \left(  \frac{j'_\textrm{1}}{{j}_\textrm{1}}  +  \frac{j'_\textrm{2}}{{j}_\textrm{2}} \right)$. However this is forbidden because for $2 \frac{j'_\textrm{av.}}{{j}_\textrm{av.}}  =  \frac{j'_\textrm{1}}{{j}_\textrm{1}}  +  \frac{j'_\textrm{2}}{{j}_\textrm{2}}$ the right hand side of 
\begin{eqnarray}
\frac{d}{dt} \left( 2 \frac{j'_\textrm{av.}}{{j}_\textrm{av.}}  -  \frac{j'_\textrm{1}}{{j}_\textrm{1}}  -  \frac{j'_\textrm{2}}{{j}_\textrm{2}}  \right)
 &=& \left(  2 \kappa_\textrm{av.} - \kappa_1 - \kappa_2 \right) - 2 \frac{(j'_\textrm{av.})^2}{{j}_\textrm{av.}^2}  +  \frac{(j'_\textrm{1})^2}{{j}_\textrm{1}^2}  +  \frac{(j'_\textrm{2})^2}{{j}_\textrm{2}^2} \\
& = &  0  - \frac{1}{2}  \left( 2 \frac{j'_\textrm{av.}}{{j}_\textrm{av.}}  -  \frac{j'_\textrm{1}}{{j}_\textrm{1}}  -  \frac{j'_\textrm{2}}{{j}_\textrm{2}}  \right) \left( 2 \frac{j'_\textrm{av.}}{{j}_\textrm{av.}}  +  \frac{j'_\textrm{1}}{{j}_\textrm{1}}  +  \frac{j'_\textrm{2}}{{j}_\textrm{2}}  \right) +    \frac{1}{2} \left(  \frac{j'_\textrm{1}}{{j}_\textrm{1}}  -    \frac{j'_\textrm{2}}{{j}_\textrm{2}} \right)^2 \nonumber 
\end{eqnarray}
is equal to $\frac{1}{2} \left(  \frac{j'_\textrm{1}}{{j}_\textrm{1}}  -    \frac{j'_\textrm{2}}{{j}_\textrm{2}} \right)^2   $ and is therefore positive.

\item Argument that it remains true that $j_\textrm{av.}^2\geq j_1  j_2$: 

if they were to cross, there would need to be a moment when $j_\textrm{av.}^2 = j_1 j_2$ and $\frac{d}{dt} j_\textrm{av.}^2 < \frac{d}{dt} j_\textrm{1} j_\textrm{2}$. However, the inequality we just established, which can be written $\frac{1}{j_\textrm{av.}^2} \frac{d}{dt} j_\textrm{av.}^2 \geq \frac{1}{j_1 j_2} \frac{d}{dt} j_\textrm{1} j_\textrm{2}$, forbids this.

\end{itemize}

     \begin{figure}[htbp] 
    \centering
    \includegraphics[width=6in]{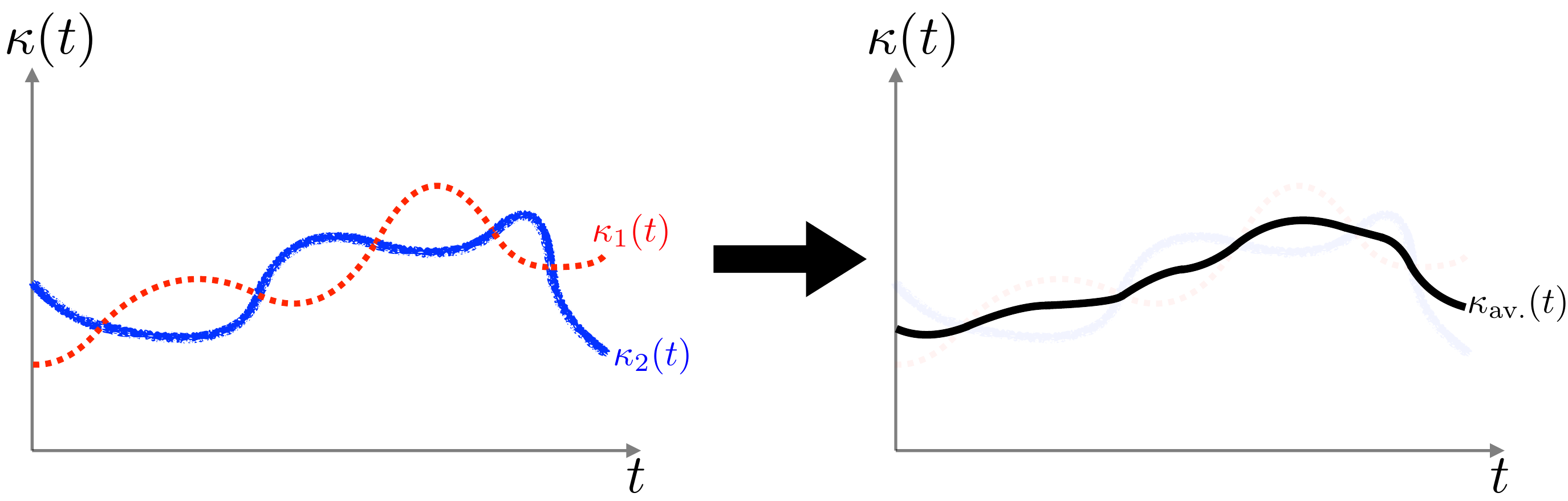} 
    \caption{The multiplicative lemma Eq.~\ref{eq:multiplicativelemma} says that averaging from $\textcolor{red}{\kappa_1(t)}$ and $\textcolor{blue}{\kappa_2(t)}$ to ${\kappa}_\textrm{av.}(t) \equiv \frac{\kappa_1(t) + \kappa_2(t)}{2}$ increases the growth of the product of the $j$s, i.e.~that $j_\textrm{av.}(t)^2 \geq \textcolor{red}{j_1(t)} \textcolor{blue}{j_2(t)}$.}
        \label{fig:kappaave}
 \end{figure}

\subsection{Averaging multiple $\kappa$s increases the product $j_1 j_2 j_3 \ldots j_{d-1} j_d$}
This follows directly from iterating the Multiplicative Lemma, Eq.~\ref{eq:multiplicativelemma}. Choose a pair of schedules $\{ \kappa_i,\kappa_j \}$ and apply Eq.~\ref{eq:multiplicativelemma}. This will increase the product $j_i j_j$ and therefore increase the product $\prod_i j_i$. Then pick a new pair $\{\kappa_k, \kappa_\ell \}$ and repeat. And then another. This process will keep increasing the product $\prod_i j_i$ and will only converge when all of the $\kappa$s are equal.

\subsection{Jacobi equation screens UV} \label{eq:secscreensUV}
We just proved that `averaging' the $\kappa$-schedules increases the product of the trajectories, but `averaging' is not an allowed primitive operation. The allowed primitive operation is `shuffling'. However, we will now show that the effect of rapidly shuffling back and forth between two $\kappa$-schedules is effectively to implement the average $\kappa$-schedule. In the language of engineering, we would say that the Jacobi equation is a low-pass filter in which high-frequency modes of $\kappa(t)$ are screened. 
Let's spell that out in detail. Imagine dividing time into $t \Delta^{-1}$ intervals each of duration $\Delta$. For odd intervals we will evolve with $\kappa_1(t)$; for even intervals we will evolve with $\kappa_2(t)$. Taylor expanding the equation of motion, and using $\kappa_2( \Delta + t) = \kappa_2(t) +  \kappa_2'(t)\Delta +   \frac{1}{2}\kappa_2''(t) \Delta^2 + \ldots$, we can write the evolution generated between $t$ and $t+2 \Delta$ as 
\begin{equation}
\left( \begin{array}{c}
j(t + 2 \Delta) \\
j'(t + 2 \Delta) \\
\end{array} \right) 
=  \left( \begin{array}{cc} 1 & 2\Delta \\ (\kappa_1 (t) + \kappa_2(t) )\Delta & 1 \end{array} \right) \left( \begin{array}{c}
j(t) \\
j'(t ) \\
\end{array} \right) 
+O(\Delta^2)  \ . 
\end{equation}
The point is that so long as $\kappa(t)$ is differentiable then the $\Delta^0$ and $\Delta^1$ terms are the same as they would have been had we instead just evolved with $\kappa_\textrm{av.}(t) = \frac{\kappa_1(t) + \kappa_2(t)}{2}$, with the first deviation being in the $\Delta^2$ term. Since as $\Delta \rightarrow 0$ the error-per-step vanishes faster than the number-of-steps increases, this proves that taking $\Delta \rightarrow 0$ gets you arbitrarily close to the evolution generated by $\kappa_{\textrm{av.}}(t)$.

\section{Our bound evaluated for example metrics} \label{appendix:examplemetrics}
Let's illustrate our new bound by evaluating it for some simple worked examples.  
\subsection{Our bound for $\mathbb{H}^2 \times \mathbb{H}^2$}  
\label{subsubsec:sanitycheckH2H2}
The example of $\mathbb{H}^2 \times \mathbb{H}^2$  is so simple that it is easy to directly calculate the volume growth, see Eq.~\ref{eq:actualforproductofhyperboles}. But for illustrative purposes we'll apply our bound, Eq.~\ref{eq:bigmainresulthere}. For  an $X^\mu$ that points a fraction $\cos \psi$ down one $\mathbb{H}^2$ and $\sin \psi$ down the other $\mathbb{H}^2$, 
\begin{eqnarray}
R &=& 1 \label{eq:Risonehere} \\
W^{a}_{\ \, b}  &=&\mathcal{R}^a_{\ \mu \nu b}X^{\mu} X^{\nu} 
+ \frac{1}{d-1}  \mathcal{R}_{\mu \nu}X^\mu X^\nu h^a_{\ b} \\
&= &  -  \left( \begin{array}{cccc}
0 &  0 & 0 & 0 \\
0 &  \cos^2 \psi  & 0 & 0 \\
0 &  0 & 0 & 0 \\
0 &  0 & 0 & \sin^2 \psi 
\end{array} \right) + \frac{1}{3}  \left( \begin{array}{cccc}
\sin^2 \psi  &  0 & - \sin \psi \cos \psi & 0 \\
0 &  1  & 0 & 0 \\
- \sin \psi \cos \psi  &  0 & \cos^2 \psi & 0 \\
0 &  0 & 0 & 1
\end{array} \right) \nonumber \\ 
\textrm{Tr}[\hat{W}^2]  &=& \frac{1}{6} + \frac{1}{2} (\cos^2 \psi - \sin^2 \psi)^2 \ . \label{eq:thevalueofW2here}
\end{eqnarray}
As required, we have $W^{a}_{\ \, b}X^b = 0$ and $\textrm{Tr}[\hat{W}] = W^{a}_{\ \, a} =0$. Because this is an Einstein metric, $R$ is the same in every direction, and because $\textrm{Tr}[\hat{W}(X)^2] \, \geq \, \frac{1}{6}$ is nonzero in every direction, our bound will give a better-than-BG bound on the growth rate.  Since this very simple example is in fact a symmetric space, $\nabla_a \mathcal{R}_{bcde} = 0$, there is no precession and 
\begin{equation}
W'_{ab} \equiv X^c \nabla_c ( \mathcal{R}_{a\mu \nu b}X^{\mu} X^{\nu} + \frac{1}{d-1}  \mathcal{R}_{\mu \nu}X^\mu X^\nu h_{ab}) = 0  \ . \label{eq:wprimeiszerohere}
\end{equation}
 Let's evaluate our bound  Eq.~\ref{eq:bigmainresulthere} for this example. Using Eqs.~\ref{eq:Risonehere}, \ref{eq:thevalueofW2here}, and \ref{eq:wprimeiszerohere}, we have 
\begin{equation}
\frac{1}{3} \langle \theta \rangle^2 \leq \max_\psi \left[ R^2 - \frac{1}{4} \frac{ \textrm{Tr}[\hat{W}^2]}{R^2}  \right]  = \max_\psi \left[ 1 - \frac{1}{4} \left( \frac{1}{6} + \frac{1}{2} (\cos^2 \psi - \sin^2 \psi)^2  \right) \right] \ . 
\end{equation} 
This is maximized in the diagonal direction $\psi = \frac{\pi}{4}$ (which is the same direction that maximizes the actual growth rate). Putting it all together, and comparing to the BG bound and the actual growth rate, Eq.~\ref{eq:actualforproductofhyperboles}, we have for this example
\begin{equation} 
 \theta_\textrm{BG}^2 = 3 > \langle \theta_\textrm{new} \rangle^2 = \frac{23}{8} >  \theta_\textrm{actual max}^2  = 2 \ . 
\end{equation}
(Since this is an Einstein metric, we could have used the enhanced bound for $R_\textrm{min} = R_\textrm{max}$ from Eq.~\ref{eq:improvementwithRmin}, which would have given the slightly tighter $\langle \theta_\textrm{Eq.~\ref{eq:improvementwithRmin}} \rangle^2 = \frac{63}{22}$.) 


\subsection{Our bound for squashed $\mathbb{H}^3$} \label{subsec:squashedH3}
The line-element of a squashed hyperbolic space is 
\begin{equation}
ds^2 = e^{2cz} dx^2 + e^{2z} dy^2 + dz^2 \ . 
\end{equation}
For $c = 1$ there is no squashing and we just have $\mathcal{H}^3$. Since $\mathcal{H}^3$ is maximally symmetric, the BG bound is exactly tight for $c=1$. For $c \neq 1$ this is a homogeneous but not isotropic space (and not a symmetric space); this is one of the left-invariant metrics on a nilpotent Lie group considered by Heintze \cite{Heintze}. Let's use polar coordinates for the unit vectors  leaving the origin
\begin{equation}
X^\mu(0) = \left( \begin{array}{c}
\cos \psi \cos \phi \\
  \cos \psi \sin \phi \\
\sin \psi  \\
\end{array} \right) \ . 
\end{equation}
The Ricci curvature in the $X^\mu$-direction is 
\begin{equation}
\mathcal{R}_{\mu \nu} X^\mu X^\nu = - c(1+c) \cos^2 \psi \cos^2 \phi - (1+ c) \cos^2 \psi \sin^2 \phi - (1+c^2) \sin^2 \psi \ . 
\end{equation}
For $c=1$ the Ricci curvature is just $-2$ in every direction. For $c>1$ the Ricci curvature is most negative in the $x$-direction, $\psi = \phi = 0$, so that $X_\textrm{min}^\mu = \{ 1, 0, 0 \}$; in that direction 
\begin{equation}
R_\textrm{max}^2 \equiv - \textrm{min}_X \mathcal{R}_{\mu \nu} X^\mu X^\nu \equiv  -  \mathcal{R}_{\mu \nu} X^\mu_\textrm{min} X^\nu_\textrm{min}=  c(1+c) \ . 
\end{equation}
Similarly, we can calculate $ \textrm{Tr}[\hat{W}^2] \equiv W_{ab}W^{ab}$, 
\begin{equation}
\textrm{Tr}[\hat{W}^2] = \frac{1}{2} (c-1)^2 \left((c+1)^2-2 (c+1) \cos ^2\psi  \left(c \sin ^2\phi +\cos ^2\phi \right)+\cos
   ^4 \psi  \left(\cos ^2\phi -c \sin ^2\phi \right)^2\right) \ . 
   \end{equation}
In the $x$-direction ($\psi = \phi = 0$) this is simply $\textrm{Tr}[\hat{W}(X_\textrm{min})^2] = \frac{1}{2} c^2 (c - 1)^2$, and the fact that $\textrm{Tr}[\hat{W}^2]$ is nonzero in the $X_\textrm{min}$ direction tells us we are guaranteed to get an improvement in the exponent over the BG bound. To evaluate how large of an improvement we will get, we need to evaluate $\textrm{Tr}[(\hat{W}')^2]$, the derivative along the geodesic, which is given by 
\begin{equation}
\textrm{Tr}[(\hat{W}')^2] = 2 (c-1)^2 c^2 \cos ^4 \psi  \left(1-\left(\sin^2 \phi- \cos ^2\phi \right)^2 \cos ^2\psi \right) \ . 
\end{equation}
This is zero in the $x$-direction, $\textrm{Tr}[(\hat{W}'(X_\textrm{min}))^2] = 0 $, and largest in the diagonal $x$-$y$-direction, $\psi = 0, \phi = \frac{\pi}{4}$, giving 
\begin{equation}
({W}_\textrm{max}')^2 \equiv \textrm{max}_{\vec{X}}\textrm{Tr}[(\hat{W}')^2] = 2 (c-1)^2 c^2 \ . 
\end{equation}
 Eq.~\ref{eq:bigmainresulthere} upperbounds the average expansion with 
\begin{eqnarray}
\frac{1}{2} \langle \theta \rangle^2 & \leq &   \textrm{max}_{\vec{X} } \Biggl[  R^2(X)  - \left(  \frac{ \sqrt{ 12 R_\textrm{max}^2  \textrm{Tr}[{\hat{W}}(X)^2]   + (W_\textrm{max}')^2  } -
   \sqrt{  (W_\textrm{max}')^2    }}{ 6 R_\textrm{max}^2  }\right)^2  \Biggl]  \ . \label{eq:improvedBGanswerinthiscase}
\end{eqnarray}
For this example it happens that the direction $\vec{X}$ that maximizes the quantity in square brackets is the same direction that maximizes $R^2 = - \mathcal{R}_{\mu \nu}X^\mu X^\nu$, namely the $x$-direction. The new bound is thus given by using $\psi = \phi = 0$, which yields
\begin{eqnarray}
\frac{1}{2} \langle \theta \rangle^2 & \leq & c(1+c) - \frac{(c-1)^2 \left(\sqrt{3 c^2+3 c+1}-1\right)^2}{18 (c+1)^2} \ . 
\end{eqnarray}
The first term is the regular BG bound, and the second term is the improvement proved in this paper. The improvement is zero for $c=1$, since in that case the space is $\mathcal{H}^3$ and the BG bound is tight. The multiplicative improvement in $\langle \theta \rangle$ gets bigger for bigger $c$, and saturates at large $c$ at $1 - \sqrt{5/6} \sim 8.7\%$. (For finite $c$, we can strengthen the improvement by considering that $R^2_\textrm{min} = 1+ c  > 0$ and using Eq.~\ref{eq:improvementwithRmin}.)

\section{Raychaudhuri worked examples} 
While the first and second Raychaudhuri equations, Eq.~\ref{eq:Raychaudhuri} and Eq.~\ref{eq:rayforshearmyversion1}, are well-known to physicists, physicists usually deploy them in their \emph{pseudo}Riemannian (spacetime) form. By contrast this paper considers purely spatial Riemannian metrics, with {all-plus} metric signature. At the same time mathematicians may be familiar with the mathematical machinery of the Raychaudhuri equation, but not the vocabulary. We'll now provide worked examples of the Raychaudhuri equation for two simple metrics.

\subsection{Example 1:   $\mathbb{H}^{d}$} \label{sec:explicitraychaudforHd}
The most trivial example is the $d$-dimensional hyperbolic space, with line element 
\begin{equation}
ds^2 = dt^2 + \sinh^2 t \, d \Omega_{d-1}^{\ 2} .
\end{equation}
The congruence of geodesics leaving the origin is purely radial, 
\begin{equation}
X_\mu  = X^\mu = \left( \begin{array}{c}
1 \\
0 \\
\vdots \\ 
0 
\end{array} \right) . 
\end{equation}
The derivative of the geodesic field is 
\begin{equation}
\nabla_\mu X_{\nu}  = \partial_{\mu} X_{\nu} - \Gamma_{\mu \nu}^{a} X_a = 0 - \Gamma^{t}_{\mu \nu} =  \Gamma^{t}_{\mu \nu}  = - \sinh t \cosh t \, h_{\mu \nu},
\end{equation}
where the projected metric is diagonal $h_{\mu}^{ \ \nu} \equiv g_{\mu}^{\  \nu} - X_\mu X^{ \nu} = \textrm{diag}\{ 0, 1, \ldots, 1 \}$. 
Let's show how this is consistent with both Raychaudhuri equations. 
\subsubsection{First Raychaudhuri: evolution of expansion for $\mathbb{H}^{d}$}
Raychaudhuri's first equation, Eq.~\ref{eq:Raychaudhuri}, is satisfied since 
\begin{eqnarray}
\theta &=&  
 (d-1) \frac{\cosh t}{\sinh t} \\ 
\theta' & = & - \frac{d-1}{\sinh^2 t} \\
\frac{\theta^2}{d-1} & = & (d-1) \frac{\cosh^2 t}{\sinh^2 t} \\
\sigma \ = \ \omega & = & 0 \\ 
\mathcal{R}_{\mu \nu} X^{\mu} X^{\nu}  \ = \  \mathcal{R}_{tt} & = & -(d-1). 
\end{eqnarray}
At $t=0$ the expansion $\theta$ is divergent, as it has to be since $\theta$ is the logarithmic derivative of the area element and the area element goes from zero to non-zero. 
At late times $\theta'$ goes to zero and the nonzero terms in the Raychaudhuri equation are $\theta^2$ and $\mathcal{R}_{tt}$. 
\subsubsection{Second Raychaudhuri: evolution of shear for $\mathbb{H}^d$}

For  $\mathbb{H}^d$, the shear of a geodesics sphere is zero. The second Raychaudhuri equation Eq.~\ref{eq:rayforshearmyversion1} is trivially satisfied because every term is zero. Shear is sourced by unequal sectional curvatures, but for  $\mathbb{H}^d$ all the sectional curvatures are the same. 

\subsection{Example 2:  $\mathbb{H}^2 \times \mathbb{H}^2$} \label{subsubsec:rayforH2H2}
The next example is the product of two unit hyperbolic planes, consider in Sec.~\ref{subsubsec:sanitycheckH2H2}, which has  line element 
\begin{equation}
ds^2 = d \tau_1^2 +  \sinh^2 {\tau_1} \, d \phi_2^2 + d \tau_2^2 + \sinh^2 {\tau_2 } \, d \phi_2^2 \ . \label{H2H2metricexplicit}
\end{equation}
Because this is an Einstein metric, the Ricci curvature is the same in every direction
\begin{equation}
\mathcal{R}_{\mu \nu} X^{\mu} X^{\nu}  = - 1 \ . 
\end{equation}
The congruence of geodesics leaving the origin is 
\begin{equation}
X^\mu = X_\mu = \frac{1}{\sqrt{\tau_1^2 + \tau_2^2}} \left( \begin{array}{c}
\tau_1 \\
0 \\
\tau_2 \\
0  
\end{array} \right) . \label{eq:H2H2pointingexplicit}
\end{equation}
The derivative of the geodesic field is 
\begin{equation}
\nabla_\mu X_{\nu}  = \partial_{\mu} X_{\nu} - \Gamma_{\mu \nu}^{a} X_a  = \partial_{\mu} X_{\nu} -  \frac{\tau_1 }{\sqrt{\tau_1^2 + \tau_2^2}} \Gamma_{\mu \nu}^{\tau_1} -  \frac{\tau_2 }{\sqrt{\tau_1^2 + \tau_2^2}} \Gamma_{\mu \nu}^{\tau_2}  \ . 
\end{equation}
Using that $\partial_{\tau_i}  X_{\tau_i} = \frac{\tau_j^2}{(\tau_1^2 + \tau_2^2)^{3/2}}$ and $\partial_{\tau_i}  X_{\tau_j} = \frac{- \tau_1 \tau_2}{(\tau_1^2 + \tau_2^2)^{3/2}}$ and $\Gamma^{\tau_i}_{\mu \nu} = - \sinh \tau_i {\cosh \tau_i} \delta_{\mu \phi_i} \delta_{\nu \phi_i}$, 
\begin{equation}
\nabla_\mu X^{\nu}  = \frac{1}{\sqrt{\tau_1^2 + \tau_2^2} }  \left( 
\begin{array}{cccc}
\frac{\tau_2^2}{\tau_1^2 + \tau_2^2} & 0 & \frac{- \tau_1 \tau_2}{\tau_1^2 + \tau_2^2} & 0 \\
0 &   \frac{\tau_1  \cosh \tau_1 }{\sinh \tau_1}   & 0 & 0 \\
\frac{- \tau_1 \tau_2}{\tau_1^2 + \tau_2^2} &0 & \frac{\tau_1^2}{\tau_1^2 + \tau_2^2} & 0 \\
0 & 0 & 0 &    \frac{\tau_2  \cosh \tau_2  }{\sinh \tau_2} \ . 
\end{array} \right) \ .   \label{eq:thederivativeofnormalforH2H2}
\end{equation}
\subsubsection{First Raychaudhuri: evolution of expansion  for $\mathbb{H}^2 \times \mathbb{H}^2$}
While the vorticity is still zero $\omega^2=0$, the shear  $\sigma^2 \equiv \textrm{Tr}[\hat{\sigma}^2] \equiv \sigma_{\mu \nu} \sigma^{\mu \nu}$ is now nonzero,
\begin{eqnarray}
\theta & = & \frac{1}{\sqrt{\tau_1^2  + \tau_2^2} } \left( 1 + \frac{\cosh \tau_1}{\sinh \tau_1} \tau_1 +\frac{\cosh \tau_2}{\sinh \tau_2} \tau_2   \right) \label{eq:hereheretheta} \\
\theta' & \equiv& X^{\mu} \nabla_\mu \theta = \frac{ \tau_1 \partial_{\tau_1} \theta + \tau_2 \partial_{\tau_2} \theta }{\sqrt{\tau_1^2 + \tau_2^2}} = -\frac{1}{\tau_1^2 + \tau_2^2} \left( 1 + \frac{\tau_1^{\, 2}}{\sinh^2 \tau_1} +\frac{\tau_2^{\, 2}}{\sinh^2 \tau_2}   \right) \\
\sigma^2 & = &  \frac{2}{3(\tau_1^2 + \tau_2^2)} \Bigl( \tau_1^2 \coth ^2 \tau_1 -\coth \tau_1 (\tau_1 \tau_2 \coth
   \tau_2+\tau_1) \nonumber \\
   &&\hspace{4.3cm}  +\tau_2 \coth \tau_2 (\tau_2 \coth \tau_2 \, -1)+1 \Bigl) \ .  \label{eq:hereheresigma} 
\end{eqnarray}
These expression satisfy Raychaudhuri's first equation, Eq.~\ref{eq:Raychaudhuri}. 
In the limit that $\tau_1$ and $\tau_2$ are very small, $\theta$ is divergent but $\sigma^2$ vanishes. 
In the limit that $\tau_1$ and $\tau_2$ get very large, $\theta'$ is the only term that does not contribute at leading order to the first Raychaudhuri equation, 
\begin{eqnarray}
 \theta & =& \frac{\tau_1 + \tau_2 }{\sqrt{\tau_1^2 + \tau_2^2} } + \ldots \\
 \theta' & = & 0 + \ldots \\ 
 \frac{1}{d-1} \theta^2 & = & \frac{1}{3} \frac{(\tau_1 + \tau_2)^2 }{{\tau_1^2 + \tau_2^2} } + \ldots \\
\sigma^2 & = & \frac{2}{3} \frac{\tau_1^2 -  \tau_1 \tau_2  +  \tau_2^2}{\tau_1^2 + \tau_2^2} + \ldots  \ . 
\end{eqnarray}
The fact that the shear is still relevant at large radius is why the Bishop-Gromov bound fails to be tight, and why there was therefore room for the improvement to the BG bound we exhibited for this metric in Sec.~\ref{subsubsec:sanitycheckH2H2}.

\subsubsection{Second Raychaudhuri: evolution of shear for $\mathbb{H}^2 \times \mathbb{H}^2$} \label{sec:shearH2H2example}
By direct calculation one can confirm that the derivative of the geodesics described in Eq.~\ref{eq:thederivativeofnormalforH2H2} satisfies the second Raychaudhuri equation, Eq.~\ref{eq:rayforshearmyversion1}. To make it compact to write down, let's instead here check that we satisfy the equation of motion for $\sigmaTWO$. This equation is derived by contracting Eq.~\ref{eq:rayforshearmyversion1} with $\sigma^{ab}$ and sum over repeated indices to give 
\begin{equation}
\frac{1}{2} \frac{ d  \sigmaTWO}{dt}  = - \frac{2}{3} \theta  \sigmaTWO - \sigmaTHREE + \sigma^{ab}  {W}_{ab}  \ . \label{eq:RAYevolutionofshear}
\end{equation}
Again 
$\sigmaTWO \equiv \sigma^2 \equiv \sigma_{a}^{\ b}\sigma_{b}^{\ a}$ and $\sigmaTHREE \equiv \sigma_{a}^{\ b}\sigma_{b}^{\ c} \sigma_{c}^{\ a}$. Let's evaluate each of the terms. 
The expansion $\theta$ is given by Eq.~\ref{eq:hereheretheta}, the shear $\sigma^2$ is given by Eq.~\ref{eq:hereheresigma}, and 
\begin{eqnarray}
   \sigmaTHREE & = & \frac{(\tau_1 \coth \tau_1-2 \tau_2 \coth \tau_2+1) (\tau_2 \coth
   \tau_2 - 2 \tau_1 \coth \tau_1 + 1) (2 - \tau_1 \coth \tau_1-\tau_2 \coth \tau_2)}{9
   \left(\tau_1^2+\tau_2^2\right)^{3/2}} \nonumber \\ 
   \sigma^{ab}  W_{ab}  
& = & \frac{1}{\left({\tau_1^2 + \tau_2^2}\right)^\frac{3}{2}} \left( \tau_1^3 \frac{ \cosh \tau_1 }{\sinh \tau_1}   + \tau_2^3 \frac{ \cosh \tau_2 }{\sinh \tau_2}   \right)  - \frac{1}{3} \theta\ . 
\end{eqnarray}
Putting these all together confirms that Eq.~\ref{eq:RAYevolutionofshear} is indeed satisfied for this example. 


\end{document}